\input ./style/arxiv-general.cfg
\documentclass[MSNbibl,number,citesort,seceqn,dvips]{arxbj}
\makeatletter
   \@ifpackageloaded{graphicx}{}{\usepackage{graphicx}}
\makeatother
\usepackage{upgreek}

\aid{0}
\volume{21}
\issue{3}
\pubyear{2015}
\firstpage{1494}
\lastpage{1537}
\doi{10.3150/14-BEJ612} % kopijuoti is 'New paper accepted'
\makeatletter

\newcommand{\rrvert}{\vert}
\newcommand{\llvert}{\vert}

\newtheorem{theorem}{Theorem}[section]
\newtheorem{lemma}{Lemma}[section]
\newremark{remark}{Remark}[section]
\newtheorem{corollary}{Corollary}[section]
\newtheorem{proposition}{Proposition}[section]
\newcommand{\eqref}[1]{(\ref{#1})}

\newcommand{\iint}{\int\!\!\!\int}

\newcommand{\Lip}{\operatorname{Lip}}
\newcommand{\tr}{\operatorname{tr}}
\newcommand{\card}{\operatorname{card}}
\newcommand{\Cov}{\operatorname{Cov}}

\renewcommand{\epsilon}{\varepsilon}
\def\sfrac#1#2{#1/#2}

\def\sklfrac#1#2{(#1/#2)}

\makeatother

\begin{document}
\begin{frontmatter}

\title{Simultaneous large deviations for the shape of Young diagrams
associated with random~words}
\runtitle{Large deviations  Young diagrams}

\begin{aug}
%%%% inicialai - be tarpu
\author{\inits{C.}\fnms{Christian}~\snm{Houdr\'e}\thanksref{e1}\ead
[label=e1,mark]{houdre@math.gatech.edu}} \and
\author{\inits{J.}\fnms{Jinyong}~\snm{Ma}\thanksref{e2}\ead
[label=e2,mark]{jma@math.gatech.edu}}
%\author{\inits{}\fnms{}~\snm{}}
%%\runauthor{} %% auto
%\dedicated{}
\address[]{Georgia Institute of Technology,
School of Mathematics, Atlanta, GA 30332-0160, USA.\\
\printead{e1,e2}}
%\address[b]{}
\end{aug}

% HISTORY:
\received{\smonth{12} \syear{2011}}
\revised{\smonth{12} \syear{2013}}

% ABSTRACT
%
\begin{abstract}
We investigate the large deviations of the shape of the random RSK
Young diagrams
associated with a random word of size $n$ whose letters are independently
drawn from an alphabet of size $m=m(n)$.
When the letters are drawn uniformly and when both $n$ and $m$ converge
together to
infinity, $m$ not growing too fast with respect to $n$, the large
deviations of
the shape of the Young diagrams are shown to be the
same as that of the spectrum of the traceless GUE. In the non-uniform
case, a control of
both highest probabilities will ensure that the length of the top row
of the diagram
satisfies a large deviation principle. In either case, both speeds and
rate functions
are identified. To complete our study, non-asymptotic concentration
bounds for the length
of the top row of
the diagrams, that is, for the length of the longest increasing
subsequence of the random word are also
given for both models.
\end{abstract}

% KEYWORDS
% visi is mazosios raides ir pagal abecele
%

\end{frontmatter}

%s1 #&#
\section{Introduction and results}

Let $\mathcal{A}_m=\{\alpha_{1}<\alpha_{2}<\cdots<\alpha_{m}\}$ be an
ordered alphabet of size $m$, and
let a word be made of the random letters $X_1^m,X_2^m,\ldots,X_n^m$,
independently drawn
from $\mathcal{A}_m$. The Robinson--Schensted--Knuth (RSK) correspondence
associates with this random word a pair of Young diagrams, of the same
shape, having at
most $m$ rows. Now for $i=1,2,\ldots, m$, let $R_i(n,m)$ denote the
length of the $i$th row of
the Young diagrams, and recall that $R_1(n,m)$, the length of
the top row, coincides with the length of the longest
increasing subsequence of the random word $X_1^mX_2^m\cdots X_n^m$.
Appropriately renormalized, and for uniform draws, the shape $
(R_i(n,m) )_{i=1}^m$ of the
Young diagrams converges, in law and with $n$, to the spectrum of an $m \times m$
element of the traceless GUE
(\cite{Jo2,TW2}). In turn, any fixed size subset of this
spectrum, also converges
with $m$, and after proper normalization, to a multidimensional Tracy--Widom
distribution (\cite{TW1,TW3}). These iterated convergence
results have further led
(see \cite{BH}) to the study of the limiting shape when both the word
length and alphabet size
simultaneously grow to infinity. This is briefly recalled next.

Let the random matrix $\mathbf{X}=(\mathbf{X}_{ij})_{1\le i,j \le m}$
be an element
of the $m \times m$ GUE with rescaling such
that $\operatorname{Re}(\mathbf{X}_{ij})\sim N(0,1/2)$ and $\operatorname{Im}(\mathbf{X}_{ij})\sim
N(0,1/2)$, for $i\ne j$;
and $\mathbf{X}_{ii}\sim N(0,1)$ (see \cite{AGZ} and \cite{Me} for
background on
random matrices). Let $(\lambda_{1}^{m},\lambda_{2}^{m},\ldots, \lambda
_{m}^{m})$ be the
non-increasing ordered spectrum of $\mathbf{X}$, and
let $(\lambda_{1}^{m,0},\lambda_{2}^{m,0},\ldots, \lambda_{m}^{m,0})$
be the corresponding
ordered spectrum of an element of the traceless GUE (i.e., of $\mathbf
{X}-\tr(\mathbf{X})/m$).
An important fact (e.g., \cite{Ba,GTW,HL2,BGH}) asserts that\vspace*{-1pt}
%
%e1.1 #&#
\begin{eqnarray}
\label{E:brownianeigen}
&& \bigl(\lambda_{1}^{m,0},\lambda_{2}^{m,0},
\ldots, \lambda _{m}^{m,0} \bigr)
\nonumber
\\[-8pt]\\[-8pt]
&&\quad \stackrel{\mathcal{L}} {=}\frac{\sqrt{m-1}}{\sqrt{m}} \boldsymbol{\Theta
}_m^{-1} \Biggl( \Biggl(\max_{\mathbf{t}\in I_{k,m}}\sum
_{j=1}^k \sum
_{l=j}^{m-k+j} \bigl(\tilde{B}_{t_{j,l}}^{l}-
\tilde {B}_{t_{j,l-1}}^{l} \bigr) \Biggr)_{1\le k \le m} \Biggr),\nonumber
\end{eqnarray}
where $\boldsymbol{\Theta}_m\dvtx  \mathbb{R}^m\to\mathbb{R}^m$ is defined via
$(\boldsymbol{\Theta}_m(\mathbf{x}))_{j}=\sum_{i=1}^j x_i, 1\le j \le m$, and
where\linebreak[4]  $(\tilde{B}^{j}_{t})_{1 \le j \le m, t \in [0,1]}$ is a
driftless $m$-dimensional Brownian motion with covariance matrix\vspace*{-1pt}
%
%e1.2 #&#
\begin{equation}
\label{covmatrix} t %
\pmatrix{ 1& \rho& \cdots& \rho
\cr
\rho& 1 & \cdots&
\rho
\cr
\vdots&\vdots&\ddots& \vdots
\cr
\rho& \rho&\cdots&1 } %
,
\end{equation}
with $\rho=-1/(m-1)$, and where for $1\le k \le m$,\vspace*{-1pt}
%
%e1.3 #&#
\begin{eqnarray*}
I_{k,m}&=&\bigl\{\mathbf{t}=(t_{j,l}\dvtx  1\le j \le k, 0
\le l\le m)\dvtx  t_{j,j-1}=0, t_{j,m-k+j}=1, 1\le j\le k,
\\
&&\hphantom{\bigl\{}t_{j,l-1}\le t_{j,l}, 1\le j\le k, 1\le l \le m-1;
t_{j,l}\le t_{j-1,l-1}, 2\le j \le k, 2 \le l \le m \bigr\}.
\end{eqnarray*}

By comparing the Brownian functionals in \eqref{E:brownianeigen} with
discrete functionals representing the shape of the random Young
diagrams, and via
a KMT approximation, under some growth conditions
on $m$, the simultaneous asymptotic convergence of
the shapes is obtained in \cite{BH}.

A related strategy is pursued here in order to investigate the large
deviations of
the shape of the RSK Young diagrams.
More precisely, we obtain a large deviation principle (LDP) for the
length of the first $r$ rows of the Young diagrams, when $n$ and $m$
simultaneously
converge to infinity and when the size $m$ of the alphabet does not
grow too fast.
To achieve our goals, we also rely on the techniques and results developed
in \cite{BDG} (see also \cite{ABC}), where large deviations are
obtained for
the largest (or the $r$th largest) eigenvalue of the GOE.
These methodologies further give the multidimensional large deviations
for the
first $r$ eigenvalues of the ordered spectrum of the traceless GUE.
In turn, when combined with a KMT approximation, these lead to large
deviations for the shape of the diagrams. Clearly, the
results presented below complement the weak convergence ones
of \cite{BH} and as any LDP results they allow to precisely quantify the
deviation from the typical (limiting) shape of Young diagrams.

Let us next put our work into context. For random permutations, the
large deviations of
the length of the longest increasing subsequence are described
in \cite{DZ} and \cite{Se}, while, moderate deviations are given
in \cite{LM} and \cite{LMR}. Closer to our framework, in \cite{Ib},
following the comparison
method of \cite{BS} and \cite{BM}, large deviations for the
last-passage directed
percolation model close to the $x$-axis are mainly
established for i.i.d. weights which are Gaussian or have
finite exponential
moments. The study of the length of the top row of the diagrams also
corresponds to a last-passage
percolation problem but with \textit{dependent} (exchangeable in the
uniform case) Bernoulli weights (see \eqref{E:tableexpress}). For
uniform draws, our framework also deals with the other
rows of the diagrams.

Here is one of our main results:
%
%th1.1 #&#
\begin{theorem}\label{T:diagramsmulti}
In the uniform case, let $m$ and $n$ simultaneously converge to
infinity in such a
way that $ m(n)=\mathrm{o}(n^{1/4})$. Then, for any $r \ge1$,\vspace*{-1pt}
\[
\biggl(\frac{R_1(n,m(n))-n/m(n)}{\sqrt{n}},\ldots,\frac
{R_r(n,m(n))-n/m(n)}{\sqrt{n}} \biggr)
\]
satisfies a large deviation principle with speed $m(n)$ and good rate
function $I_r$ on the
space $\mathcal{L}^r:=\{(x_1,x_2,\ldots,x_r)\in\mathbb{R}^r\dvtx x_1\ge x_2
\ge\cdots\ge x_r\}$, where\vspace*{-1pt}
%
%e1.4 #&#
\begin{equation}
\label{E:ratefunction} I_r(x_1,x_2,
\ldots,x_r)= %
\cases{ 2\displaystyle \sum_{i=1}^r
\int_2^{x_i}\sqrt{(z/2)^2-1}\,
\mathrm{d}z, &\quad \mbox{if $x_1\ge x_2\ge \cdots\ge
x_r \ge2$,}
\cr
+\infty, &\quad \mbox{otherwise.}} %
\end{equation}
In other words, for all $x_1\ge x_2\ge \cdots\ge x_r \ge2$,\vspace*{-1pt}
%
%e1.5 #&#
%e1.6 #&#
\begin{eqnarray}\label{tablemulti:finite}
&&\lim_{n \to\infty} \frac{1}{m(n)} \log\mathbb{P} \biggl(
\frac
{R_1(n,m(n))-n/m(n)}{\sqrt{n}} \ge x_1,\ldots,\frac{R_r(n,m(n))-n/m(n)}{\sqrt{n}} \ge x_r \biggr)
\nonumber\\[-8pt]\\[-8pt]
&&\quad =-2\sum
_{i=1}^r \int_2^{x_i}
\sqrt{(z/2)^2-1}\,\mathrm{d}z,\nonumber
\end{eqnarray}
while for any $x<2$ and $1\le i \le r$,\vspace*{-1pt}
%
%e1.7 #&#
\begin{equation}
\label{tablemulti:infinity} \lim_{n \to\infty}\frac{1}{m(n)}\log\mathbb{P}
\biggl(\frac{R_i(n,m(n))-n/m(n)} {
\sqrt{n}}\le x \biggr)=-\infty.
\end{equation}
\end{theorem}
%
%re1.1 #&#
\begin{remark}\label{R:1}
\textup{(i)} The rate function $I_r$ in \eqref{E:ratefunction} is a good rate function.
Moreover, it is continuous and increasing with respect to each
individual variable on its effective
domain $\mathcal{D}_{I_r}=\{(x_1,x_2,\ldots,x_r)\in\mathbb{R}^r\dvtx x_1\ge
x_2 \ge\cdots\ge x_r \ge2\}$,
given that the other variables are fixed. Therefore, to prove a large
deviation principle (LDP) as
in Theorem~\ref{T:diagramsmulti}, it is enough to prove a limiting
equality on rectangular subsets as in \eqref{tablemulti:finite} or
\eqref{tablemulti:infinity}
instead of proving both the usual upper and lower bounds, that
is, that
for any closed set $F$ in $\mathcal{L}^r=\{(x_1,x_2,\ldots,x_r)\in
\mathbb{R}^r\dvtx x_1\ge x_2
\ge\cdots\ge x_r\}$,\vspace*{-1pt}
%
%e1.8 #&#
\begin{equation}
\label{E:LDPF} \limsup_{n \to\infty}\frac{1}{m(n)}\log\mathbb{P}
\bigl(X_r^{n} \in F \bigr)\le-\inf_{F}I_r,
\end{equation}
and that for any open set $O$ in $\mathcal{L}^r$,\vspace*{-1pt}
%
%e1.9 #&#
\begin{equation}
\label{E:LDPO} \liminf_{n \to\infty}\frac{1}{m(n)}\log\mathbb{P}
\bigl(X_r^{n} \in O \bigr)\ge-\inf_{O}I_r,
\end{equation}
where\vspace*{-1pt}
\[
X_r^{n}= \biggl(\frac{R_i(n,m(n))-n/m(n)}{\sqrt{n}} \biggr)_{1\le i \le r}.
\]

\textup{(ii)} The restriction $m=\mathrm{o}(n^{1/4})$ (or $m=\mathrm{o}(n^{1/6})$ below) is a
technical one
and there is no reason to believe it is sharp. One can envision that
our results
continue to hold under at least a condition such as $m=\mathrm{o}(\sqrt n)$.
\end{remark}

In Theorem~\ref{T:diagramsmulti}, if at least one of the renormalized
variables is on the left
of its simultaneous asymptotic mean, by changing the convergence speed
from $m$ to $m^2$, a
more accurate form of (\ref{tablemulti:infinity}) is valid.
Below, the closed form expression obtained for $K$
was found after Satya Majumdar kindly suggested that the methodology developed
in \cite{NMV} would apply to our traceless GUE framework.
%
%th1.2 #&#
\begin{theorem}\label{P:tableautraceless}
In the uniform case, let $m$ and $n$ simultaneously converge to
infinity in such
a way that $m(n)=\mathrm{o}(n^{1/6})$. Then, for any $r \ge1$,
\[
\biggl(\frac{R_1(n,m(n))-n/m(n)}{\sqrt{n}},\ldots,\frac
{R_r(n,m(n))-n/m(n)}{\sqrt{n}} \biggr)
\]
satisfies a large deviation principle with speed $(m(n))^2$ and good
rate function $K(x_r)$ on the space $\mathcal{L}^r:=\{(x_1,x_2,\ldots
,x_r)\in\mathbb{R}^r\dvtx x_1\ge x_2 \ge\cdots\ge x_r\}$, where $K$ is
the rate function in the large deviation principle for
the largest eigenvalue of the $m \times m$ traceless GUE, when
on the left of its asymptotic mean. It is given by
%
%e1.10 #&#
\begin{equation}
\label{ratefunction_traceless} K(x):=\inf_{\mu\in\mathcal{M}_{0}((-\infty,x])}I(\mu),
\end{equation}
where $I$ (see \textup{\eqref{traceless:ratefunction}}) is the rate function in
the large deviation principle for
the spectral measure of the GUE, and $\mathcal{M}_{0}((-\infty,x])$ is
the set of
zero mean probability measures supported on $(-\infty,x]$.
For $x\le0$, $K(x)=+\infty$, for $x\ge2$, $K(x)=0$, and for $0<x<2$,
%
%e1.11 #&#
%e1.12 #&#
%e1.13 #&#
%e1.14 #&#
\begin{eqnarray}
\label{ratefunction_explicit} K(x)&=&\frac{1}{48} \bigl(3 \bigl(9 \sqrt[3]{2}
3^{2/3} \bigl(\sqrt{81 x^2+12}-9 x \bigr)^{2/3}-8
\bigr) x^2\nonumber
\\
&&\hphantom{\frac{1}{48} \bigl(}{}+9 \sqrt[3]{2} \sqrt[6]{3} \bigl(\sqrt{81 x^2+12}-9 x
\bigr)^{1/3} \nonumber\\
&&\hphantom{\frac{1}{48} \bigl(+9}{}\times\bigl(\sqrt{27 x^2+4} \bigl(\sqrt{81
x^2+12}-9 x \bigr)^{1/3}-5 \sqrt[3]{2} \sqrt[6]{3} \bigr) x\nonumber
\\
&&\hphantom{\frac{1}{48} \bigl(}{}-6\sqrt[3]{2}\, 3^{2/3} \bigl(\sqrt{81 x^2+12}-9 x
\bigr)^{2/3}\\
&&\hphantom{\frac{1}{48} \bigl(}{}-3\, 2^{2/3} 3^{5/6} \sqrt{27
x^2+4} \bigl(\sqrt{81 x^2+12}-9 x \bigr)^{1/3}\nonumber
\\
&&\hphantom{\frac{1}{48} \bigl(}{}+16 \log \bigl( \sqrt{81 x^2+12}-9 x \bigr)\nonumber\\
&&\hphantom{\frac{1}{48} \bigl(}{}-48 \log \bigl( 2\sqrt
[3]{3}-\sqrt[3]{2} \bigl( \sqrt{81 x^2+12}-9 x \bigr)^{2/3}
\bigr)+ 60+32 \log6 \bigr).\nonumber
\end{eqnarray}

In other words, for all $x_r\le x_{r-1} \le\cdots\le x_1 $, with
$x_r\le2$,
%
%e1.15 #&#
%e1.16 #&#
\begin{eqnarray}
\label{th2_K} &&\lim_{n \to\infty} \frac{1}{(m(n))^2} \log\mathbb{P}
\biggl(\frac
{R_1(n,m(n))-n/m(n)}{\sqrt{n}} \le x_1,\ldots,
\nonumber\\[-8pt]\\[-8pt]
&&\hphantom{\lim_{n \to\infty} \frac{1}{(m(n))^2} \log\mathbb{P}
\biggl(}\frac{R_r(n,m(n))-n/m(n)}{\sqrt{n}} \le x_r \biggr) =-K(x_r),\nonumber
\end{eqnarray}
while for all $2\le x_r\le x_{r-1} \le\cdots\le x_1$,
%
%e1.17 #&#
%e1.18 #&#
\begin{eqnarray}
\label{th2_zero} &&\lim_{n \to\infty} \frac{1}{(m(n))^2} \log\mathbb{P}
\biggl(\frac
{R_1(n,m(n))-n/m(n)}{\sqrt{n}} \le x_1,\ldots,
\nonumber\\[-8pt]\\[-8pt]
&&\hphantom{\lim_{n \to\infty} \frac{1}{(m(n))^2} \log\mathbb{P}
\biggl(}\frac{R_r(n,m(n))-n/m(n)}{\sqrt{n}} \le x_r \biggr) =0.\nonumber
\end{eqnarray}
\end{theorem}

The LDP for the longest increasing subsequence is now a simple consequence.
%
%co1.1 #&#
\begin{corollary}\label{coro3}
Let $m$ and $n$ simultaneously converge to infinity in such a way that
$m(n)=\mathrm{o}(n^{1/4})$, then
for any $x \ge2$,
\[
\lim_{n \to\infty} \frac{1}{m(n)} \log\mathbb{P} \biggl(
\frac
{R_1(n,m(n))-n/m(n)}{\sqrt{n}} \ge x \biggr)=-2 \int_2^{x}
\sqrt{(z/2)^2-1}\,\mathrm{d}z,
\]
and similarly, if $m(n)=\mathrm{o}(n^{1/6})$, for any $x \le2$,
\[
\lim_{n \to\infty} \frac{1}{(m(n))^2} \log\mathbb{P} \biggl(
\frac
{R_1(n,m(n))-n/m(n)}{\sqrt{n}} \le x \biggr)=-K(x).
\]
\end{corollary}
%
%re1.2 #&#
\begin{remark}
The methodologies developed in this paper allow to derive LDPs in
related problems.
Such is the case for last-passage directed percolation close to the
$x$-axis, or for the
departure time from many queues in series when the number of customers
is a fractional power of
the number of servers. In these two problems, similar (discrete)
functional representations
are available but with i.i.d. weights, so the large deviations rate
functions should be
the corresponding rate functions in the LDP for the largest eigenvalue
of the GUE.
\end{remark}

When the independent random letters are no longer uniformly drawn,
let $X_i^m, 1\le i \le n$, be independently and identically distributed
with $\mathbb{P}(X_1^m=\alpha_j)=p_j^m$, $1\le j \le m$.
Moreover, let $p_{\mathrm{max}}^m=\max_{1 \le j \le m} p_j^m $,
let $p_{\mathrm{2nd}}^m=\max\{p_j^m<p_{\mathrm{max}}^m\dvtx  1\le j \le m\}$, and let
also $J(m)=\{j\dvtx p_j^m=p_{\mathrm{max}}^m\}$, where $k(m)=\card (J(m))$, that is, $k(m)$
is the multiplicity of $p_{\mathrm{max}}^m$.
%
%th1.3 #&#
\begin{theorem}\label{T:nonuniform}
In the non-uniform case, let $k(m(n))$ and $n$ simultaneously converge
to infinity in
such a way that $k(m(n))^3/p_{\mathrm{max}}^m=\mathrm{o}(n)$, and let
%
%e1.19 #&#
\begin{equation}
\label{E:intheorem} \frac{n(p_{\mathrm{2nd}}^m)^2}{k(m(n))p_{\mathrm{max}}^m}=\mathrm{o}\bigl(\exp\bigl(-k\bigl(m(n)
\bigr)^{\alpha}\bigr)\bigr),\qquad  \mbox{for some } \alpha>1,
\end{equation}
where $p_{\mathrm{2nd}}$ is the second highest probability. Then,
\[
\frac{R_1(n,m(n))-np_{\mathrm{max}}^m}{\sqrt{nk(m(n))p_{\mathrm{max}}^m}}
\]
satisfies a LDP on $\mathbb{R}$ with speed $k(m(n))$ and good rate
function $I_1$.

In other words, for any $x\ge2$,
%
%e1.20 #&#
\begin{equation}
\label{E:nonuniformresult1} \lim_{n \to\infty}\frac{1}{k(m(n))} \log\mathbb{P}
\biggl(\frac
{R_1(n,m(n))-np_{\mathrm{max}}^m}{\sqrt{nk(m(n))p_{\mathrm{max}}^m}} \ge x \biggr)=-2 \int_2^{x}
\sqrt{(z/2)^2-1}\,\mathrm{d}z,
\end{equation}
while for any $x<2$,
%
%e1.21 #&#
\begin{equation}
\label{E:nonuniformresult0} \lim_{n \to\infty}\frac{1}{k(m(n))} \log\mathbb{P}
\biggl(\frac
{R_1(n,m(n))-np_{\mathrm{max}}^m} {
\sqrt{nk(m(n))p_{\mathrm{max}}^m}} \le x \biggr)=-\infty.
\end{equation}
\end{theorem}
%
%re1.3 #&#
\begin{remark}
\textup{(i)} Above, the condition $k(m(n))^3/p_{\mathrm{max}}^m=\mathrm{o}(n)$ matches exactly the
condition
$m=\mathrm{o}(n^{1/4})$ of Theorem~\ref{T:diagramsmulti}, but the new
condition \eqref{E:intheorem} is not present there. A similar remark applies
to Theorem~\ref{P:tableautraceless} and to Theorem~\ref{P:nonuniform}.

\textup{(ii)} In contrast to our first theorem, the one just stated is only for
the first row of the diagrams and not for the whole shape.
For non-uniform draws, a LDP shape result is also possible under conditions
involving all the distinct probabilities and their respective multiplicity.
These conditions are rather involved
and the corresponding proofs rather
tedious; therefore only a first row result is given above.
A similar remark applies to Theorem~\ref{P:tableautraceless} and to
Theorem~\ref{P:nonuniform}.
\end{remark}

When the renormalized variable is on the left of its simultaneous
asymptotic mean, again
a more accurate form of \eqref{E:nonuniformresult0} is possible.
Before presenting this statement, let us first recall a few facts.
For the alphabet $\mathcal{A}_m$ with corresponding set of probabilities
$\mathcal{P}=\{p_1^m,p_2^m,\ldots,p_m^m\}$,
let $p^{(1)}>p^{(2)}>\cdots>p^{(l)}$, $1\le l \le m$, be the distinct elements
in $\mathcal{P}$, and let $d_1,\ldots,d_l$ be the corresponding
multiplicities, with $\sum_{i=1}^l d_i=m$.
Then $p^{(1)}=p_{\mathrm{max}}^m$ and $d_1=k(m)$ as in the previous notations.
Let $\mathcal{G}_{m}(d_1,\ldots,d_l)$ be the set of
$m\times m$ random matrices $\mathbf{X}$ which are direct sums of mutually
independent elements of the $d_i \times d_i$ GUE, $1\le i\le l$.
Moreover, let $p_{(1)}\ge p_{(2)}\ge\cdots\ge p_{(m)}$ be the non-increasing
rearrangement of $\mathcal{P}$. The ``generalized'' $m \times m$
traceless GUE associated
with $\mathcal{P}$ is the set, denoted
by $\mathcal{G}^0(p_1^m,p_2^m,\ldots,p_m^m)$, of $m \times m$ random matrices
$\mathbf{X^0}$, of the form
%
%e1.22 #&#
\begin{equation}
\label{gene} \mathbf{X}_{i,j}^0= %
\cases{
\mathbf{X}_{i,i}-\sqrt{p_{(i)}}\displaystyle \sum
_{h=1}^m \sqrt{p_{(h)}}\,
\mathbf{X}_{h,h},  &\quad\mbox{if $i=j$,}
\cr
\mathbf{X}_{i,j}, &\quad
\mbox{otherwise,}} %
\end{equation}
where $\mathbf{X} \in\mathcal{G}_{m}(d_1,\ldots,d_l)$. Finally,
let $\tilde{\lambda}_1^0$ be the
largest eigenvalue of the diagonal block corresponding
to $p^{(1)}=p_{\mathrm{max}}^m$ in $\mathbf{X^0}$.
%
%th1.4 #&#
\begin{theorem}\label{P:nonuniform}
Let $k(m(n))$ and $n$ simultaneously converge to infinity in such a way
that $k(m(n))^5/p_{\mathrm{max}}^m=\mathrm{o}(n)$, let also
%
%e1.23 #&#
\begin{equation}
\label{E:intheorem1} \frac{n(p_{\mathrm{2nd}}^m)^2}{k(m(n))p_{\mathrm{max}}^m}=\mathrm{o}\bigl(\exp\bigl(-k\bigl(m(n)
\bigr)^{\alpha}\bigr)\bigr),\qquad  \mbox{for some } \alpha>2,
\end{equation}
where $p_{\mathrm{2nd}}$ is the second highest probability,
and let,
%
%e1.24 #&#
\begin{equation}
\label{E:intheorem2} \lim_{n \to\infty} k\bigl(m(n)\bigr)p_{\mathrm{max}}^m
=\eta,
\end{equation}
for some $0\le\eta\le1$. Then,
\[
\frac{R_1(n,m(n))-np_{\mathrm{max}}^m}{\sqrt{nk(m(n))p_{\mathrm{max}}^m}}
\]
satisfies a LDP on $\mathbb{R}$ with speed $(k(m(n)))^2$ and good rate function
$K_{\eta}$, where $K_{\eta}$ is the rate function of $\tilde{\lambda
}_1^0$ when on
the left of its asymptotic mean.

In other words, for any $ x\le2$,
%
%e1.25 #&#
\begin{equation}
\label{E:nonuniformresult2} \lim_{n \to\infty}\frac{1}{(k(m(n)))^2} \log\mathbb{P}
\biggl(\frac{R_1(n,m(n))-np_{\mathrm{max}}^m}{\sqrt{nk(m(n))p_{\mathrm{max}}^m}} \le x \biggr)=-K_{\eta}(x),
\end{equation}
while for any $x \ge2$,
%
%e1.26 #&#
\begin{equation}
\label{E:nonuniformresult222} \lim_{n \to\infty}\frac{1}{(k(m(n)))^2} \log\mathbb{P}
\biggl(\frac{R_1(n,m(n))-np_{\mathrm{max}}^m}{\sqrt {nk(m(n))p_{\mathrm{max}}^m}} \le x \biggr)=0.
\end{equation}
\end{theorem}
%
%re1.4 #&#
\begin{remark}\label{remark1}
The rate function $K_{\eta}$ is given by
\[
K_{\eta}(x)=\sup_{y\le0} \biggl(xy-yS(y)+J\bigl(S(y)
\bigr)+\frac{\eta y^2}{2} \biggr),
\]
where $J$ is the rate function (with speed $m^2$) of the largest
eigenvalue of the $m \times m$ GUE, and for each $y\le0$, $S(y)$ is
the unique
solution to $J^\prime(t)=y$ with $t\le2$.
For $x \ge2$, $J(x)=0$, while for $x\le2$, the following closed form
expression for $J$ is obtained in \cite{DM},
%
%e1.27 #&#
%e1.28 #&#
\begin{eqnarray}
\label{explicitforJ} J(x)&=&\frac{1}{216} \biggl(-x \bigl(-72x+x^3+30
\sqrt{12+x^2}+x^2\sqrt {12+x^2} \bigr)
\nonumber\\[-8pt]\\[-8pt]
&&\hphantom{\frac{1}{216} \biggl(}{}-216\log \biggl(\frac{1}{6} \bigl(x+\sqrt{12+x^2} \bigr)
\biggr) \biggr).\nonumber
\end{eqnarray}
In particular, $K_0=J$ and $K_1=K$.
In fact the relationship between the spectrum of GUE and traceless GUE
implies that
\[
K(x)=\sup_{y\le0} \biggl(xy-J^*(y)+\frac{y^2}{2} \biggr),
\]
where $*$ denotes the Legendre transform.
For any $0\le\eta\le1$, $K_{\eta}(x)=0$, when $x\ge2$, while for
$0\le\eta<1$
and $x\in(-\infty,2)$, $K_{\eta}(x)>0$ and is asymptotically
equivalent to
\[
\frac{x^2}{2 (1-\eta)}+\log \biggl(-\frac{x}{1-\eta} \biggr),
\]
as $x\to-\infty$. For $\eta=1$, when $0<x<2$, $K_1(x)=K(x)$ is
positive and finite.
As $x\to0$, $K(x)\sim-\log x$, while as $x\to2$, $K(x)\sim
C(2-x)^3$, for some
positive constant $C$.
\end{remark}

To complement the previous results, we provide corresponding
concentration results.
These rely in part on deviations inequalities for the
largest eigenvalue of the $m\times m$ GUE matrix, obtained respectively,
in \cite{Au} and \cite{LR}.
Comparing the forthcoming result with Corollary~\ref{coro3}, we see
that, in this case,
the deviation rates match the fluctuation results.
In turn these rates match the order of the tails of the
Tracy--Widom distribution. Our results are only stated for the top row
of the
diagrams, but the techniques easily apply to the whole shape when combined
with deviation inequalities for the whole spectrum of the GUE.
%
%th1.5 #&#
\begin{theorem}\label{T:concentration}
In the uniform model, let $0<\alpha< 1/4$, and
let $m\le An^{\alpha}$, for some $A>0$.
Then, for any $0<\epsilon<1$,
%
%e1.29 #&#
\begin{equation}
\label{E:concentrationthm1} \mathbb{P} \biggl(\frac{R_1(n,m)-n/m}{\sqrt{n/m}} \ge2\sqrt{m}(1+\epsilon )
\biggr) \le C(A,\alpha) \exp \biggl(-\frac{m\epsilon^{3/2}}{C(A,\alpha)} \biggr),
\end{equation}
where
\[
C(A,\alpha)= C\max\bigl(A^{10/3},1\bigr) \frac{1+\alpha}{1-4\alpha} \exp \biggl(
\frac{1+\alpha}{1-4\alpha} \biggr),
\]
for some absolute constant $C>0$.

Likewise, let $0<\alpha< 1/6$, and let $m\le An^{\alpha}$, for some $A>0$.
Then, for any $0<\epsilon<1$,
%
%e1.30 #&#
\begin{equation}
\label{E:concentrationthm2} \mathbb{P} \biggl(\frac{R_1(n,m)-n/m}{\sqrt{n/m}} \le2\sqrt{m}(1-\epsilon )
\biggr)\le{\bar C}(A,\alpha) \exp \biggl(-\frac{m^2\epsilon^{3}}{{\bar C}(A,\alpha)} \biggr),
\end{equation}
where
\[
{\bar C}(A,\alpha)=C \max\bigl(A^4,1\bigr)\frac{1+\alpha}{1-6\alpha} \exp
\biggl(\frac{1+\alpha}{1-6\alpha} \biggr),
\]
for some absolute constant $C>0$.
\end{theorem}

Again, in the non-uniform case, similar results hold
but under a further control of the second highest probability.
%
%th1.6 #&#
\begin{theorem}\label{T:nonuniformconcentration}
In the non-uniform model, let $\alpha>3$, let $k(m(n))^{\alpha
}/p_{\mathrm{max}}^m\le An$,
for some $A>0$ and let
%
%e1.31 #&#
\begin{equation}
\label{E:nonconexp1} \frac{n(p_{\mathrm{2nd}}^m)^2}{k(m(n))p_{\mathrm{max}}^m}\le B \exp\bigl(-k\bigl(m(n)\bigr)\bigr),
\end{equation}
for some $B>0$. Then, for any $0<\epsilon<1$,
%
%e1.32 #&#
\begin{equation}
\label{E:nonuniformconcentrationthm1} \mathbb{P} \biggl(\frac{R_1(n,m)-np_{\mathrm{max}}^m}{\sqrt{nk(m)p_{\mathrm{max}}^m}} \ge 2(1+\epsilon) \biggr)\le
C(A,B,\alpha) \exp \biggl(-\frac{k(m)\epsilon
^{3/2}}{C(A,B,\alpha)} \biggr),
\end{equation}
where
\[
C(A,B,\alpha)=C \max\bigl(A^{10/3\alpha},1\bigr)\max(\sqrt{B},1)
\frac{\alpha
+2}{\alpha-3} \exp \biggl(\frac{\alpha+2}{\alpha-3} \biggr),
\]
for some absolute constant $C>0$.

Likewise, let $\alpha>5$, let $k(m(n))^{\alpha}/p_{\mathrm{max}}^m\le An$, for
some $A>0$, and let
%
%e1.33 #&#
\begin{equation}
\label{E:nonconexp2} \frac{n(p_{\mathrm{2nd}}^m)^2}{k(m(n))p_{\mathrm{max}}^m}\le B\exp\bigl(-k\bigl(m(n)\bigr)^2
\bigr),
\end{equation}
for some $B>0$. Then, for any $0<\epsilon<1$,
%
%e1.34 #&#
\begin{equation}
\label{E:nonuniformconcentrationthm2} \mathbb{P} \biggl(\frac{R_1(n,m)-np_{\mathrm{max}}^m}{\sqrt{nk(m)p_{\mathrm{max}}^m}} \le2(1-\epsilon) \biggr)\le{
\bar C}(A,B,\alpha) \exp \biggl(-\frac
{k(m)^2\epsilon^{3}}{{\bar C}(A,B,\alpha)} \biggr),
\end{equation}
where
\[
{\bar C}(A,B,\alpha)=C \max\bigl(A^{4/\alpha},1\bigr) \max(\sqrt{B},1)
\frac{\alpha+2}{\alpha-5} \exp \biggl(\frac{\alpha
+2}{\alpha-5} \biggr),
\]
for some absolute constant $C>0$.
\end{theorem}

%s2 #&#
\section{Proof of Theorem \texorpdfstring{\protect\ref{T:diagramsmulti}}{1.1} 
and Theorem \texorpdfstring{\protect\ref
{P:tableautraceless}}{1.2}}
As in \cite{BH}, let
%
%e2.1 #&#
\begin{equation}
X_{i,j}^m= %
\cases{ 1, &\quad \mbox{if
$X_i^m=\alpha_j$,}
\cr
0, &\quad \mbox{otherwise,} } %
\end{equation}
be Bernoulli random variables with parameter $1/m$.
For a fixed $j$, $1 \le j \le m$, the $X_{i,j}^m$'s are i.i.d. while for
$j \neq j^\prime$, $(X_{1,j}^m,\ldots,X_{n,j}^m)$ and
$(X_{1,j^\prime}^m,\ldots,X_{n,j^\prime}^m)$ are identically distributed
but no longer independent.

Let $S_k^{m,j}=\sum_{i=1}^k X_{i,j}^m$ be the number of occurrences of
$\alpha_j$
among $(X_i^m)_{1\le i\le k}$. Since for $1\le k<l \le n$, the number of
occurrences of $\alpha_j$ among $(X_i^m)_{k+1\le i\le l}$ is
$S_l^{m,j}-S_k^{m,j}$,
\[
R_{1}(n,m)=\sup_{0=l_{0} \le l_{1}
\le\cdots \le l_{m}=n} \sum
_{j=1}^{m}\bigl(S_{l_j}^{m,j}-S_{l_{j-1}}^{m,j}
\bigr),
\]
with the convention that $S_0^{m,j}=0$.

Moreover, letting $V_k(n,m)=\sum_{i=1}^k R_i(n,m)$, combinatorial
arguments yield (see Theorem~3.1 in \cite{HL2})
%
%e2.2 #&#
\begin{equation}
\label{def} V_{k}(n,m)=\sup_{\mathbf{t}\in I_{k,m}(n)}\sum
_{j=1}^k \sum_{l=j}^{m-k+j}
\bigl(S_{[t_{j,l}]}^{m,l}-S_{[t_{j,l-1}]}^{m,l} \bigr), \qquad 1
\le k \le m,
\end{equation}
where
%
%e2.3 #&#
%e2.4 #&#
\begin{eqnarray*}
I_{k,m}(n)&=&\bigl\{\mathbf{t}=(t_{j,l}\dvtx  1\le j \le
k, 0\le l\le m)\dvtx  t_{j,j-1}=0, t_{j,m-k+j}=n,
\\
&&\hphantom{\bigl\{}1\le j\le k; t_{j,l-1}\le t_{j,l}, 1\le j\le k, 1\le l \le
m-1;
\\
&&\hphantom{\bigl\{}t_{j,l}\le t_{j-1,l-1}, 2\le j \le k, 2 \le l \le m \bigr\}.
\end{eqnarray*}
Let $\tilde{X}_{i,j}^m=(X_{i,j}^m-1/m)/\sigma_m$, with $\sigma
_m^2=(1/m)(1-1/m)$,
let $\tilde{S}_k^{m,j}=\sum_{i=1}^k \tilde{X}_{i,j}^m$.
Similarly define $\tilde{V}_{k}(n,m), 1\le k \le m$ and let $\tilde
{R}_{k}(n,m)=\tilde{V}_{k}(n,m)-\tilde{V}_{k-1}(n,m)$, $2\le k \le m$, while
$\tilde{R}_{1}(n,m)=\tilde{V}_{1}(n,m)$. Clearly $V_{k}(n,m)=\sigma
_m\tilde{V}_{k}(n,m)+kn/m$, and
%
%e2.5 #&#
\[
\frac{R_k(n,m)-n/m}{\sqrt{n}}=\sqrt{1-\frac{1}{m}}\frac{\tilde
{R}_k(n,m)}{\sqrt{nm}}.
\]
Let
%
%e2.6 #&#
\begin{equation}
\label{E:tableexpress} \tilde{V}_{k}(n,m)=\sup_{\mathbf{t}\in I_{k,m}(n)}\sum
_{j=1}^k \sum
_{l=j}^{m-k+j} \bigl(\tilde{S}_{[t_{j,l}]}^{m,l}-
\tilde{S}_{[t_{j,l-1}]}^{m,l} \bigr), \qquad 1\le k \le m,
\end{equation}
with
%
%e2.7 #&#
\begin{equation}
\Cov\bigl(\tilde{S}_{\ell}^{m,i},\tilde{S}_{\ell}^{m,j}
\bigr)= %
\cases{ \ell, &\quad \mbox{if $i=j$},
\cr
\rho\ell, &\quad \mbox{otherwise,} } %
\end{equation}
and $\rho=-1/(m-1)$.

Next, $\tilde{V}_{k}(n,m)$ can be approximated by
%
%e2.8 #&#
\begin{equation}
\tilde{L}_{k}(n,m)=\sup_{\mathbf{t}\in I_{k,m}(n)}\sum
_{j=1}^k \sum_{l=j}^{m-k+j}
\bigl(\tilde{B}_{t_{j,l}}^{l}-\tilde {B}_{t_{j,l-1}}^{l}
\bigr), \qquad 1\le k \le m,
\end{equation}
where $(\tilde{B}^{j})_{1 \le j \le m}$ is a driftless $m$-dimensional
Brownian motion with covariance matrix given in \eqref{covmatrix}, and
\[
\tilde{L}_k(n,m)\stackrel{\mathcal{L}} {=}\sqrt{n}
\tilde{L}_k(1,m).
\]

More precisely, inspired by \cite{BM},
%
%e2.9 #&#
\begin{equation}
\label{E} \bigl|\tilde{V}_{k}(n,m)-\tilde{L}_{k}(n,m)\bigr|\le2k
\sum_{l=1}^m\bigl(Y_n^{m,l}+W_n^l
\bigr),
\end{equation}
where
\[
Y_n^{m,l}=\max_{1\le i \le n}\bigl|
\tilde{S}_{i}^{m,l}-\tilde{B}_{i}^{l}\bigr|
\quad \mbox{and}\quad  W_n^l=\mathop{\sup_{ 0\le s,t \le n}}_{ |s-t|\le1}
\bigl|\tilde{B}_{s}^{l}-\tilde{B}_{t}^{l}\bigr|.
\]
Since
\[
\bigl(\tilde{R}_{k}(n,m) \bigr)_{1\le k \le m}=\boldsymbol{\Theta
}_m^{-1} \bigl( \bigl(\tilde{V}_{k}(n,m)
\bigr)_{1\le k \le m} \bigr),
\]
for any $\epsilon>0$, and from \eqref{E},
%
%e2.10 #&#
\begin{eqnarray}\label{E:concentrationback}
&&\mathbb{P} \bigl(\bigl\llvert \tilde{R}_{k}(n,m)- \bigl(
\tilde{L}_k(n,m)-\tilde{L}_{k-1}(n,m) \bigr)\bigr\rrvert \ge
\sqrt {mn}\epsilon \bigr)
\nonumber
\\
&&\quad \le\mathbb{P} \Biggl(2(2k-1)\sum_{l=1}^m
\bigl(Y_n^{m,l}+W_n^l\bigr)\ge\sqrt
{mn}\epsilon \Biggr)
\nonumber
\\
&&\quad \le\mathbb{P} \Biggl(\sum_{l=1}^mY_n^{m,l}
\ge \frac{\sqrt{mn}\epsilon}{4(2k-1)} \Biggr)+\mathbb{P} \Biggl(\sum
_{l=1}^mW_n^l\ge
\frac{\sqrt{mn}\epsilon}{4(2k-1)} \Biggr)
\\
&&\quad \le\sum_{l=1}^m \biggl(\mathbb{P}
\biggl(Y_n^{m,l}\ge\frac{\sqrt {mn}\epsilon}{m(8k-4)} \biggr)+\mathbb{P}
\biggl(W_n^{l}\ge\frac{\sqrt {mn}\epsilon}{m(8k-4)} \biggr) \biggr)
\nonumber
\\
&&\quad =m\mathbb{P} \biggl(Y_n^{m,1}\ge\frac{\sqrt{n}\epsilon}{\sqrt {m}(8k-4)}
\biggr)+m\mathbb{P} \biggl(W_n^{1}\ge\frac{\sqrt{n}\epsilon
}{\sqrt{m}(8k-4)}
\biggr), \nonumber
\end{eqnarray}
for $1\le k \le m$, and with the convention
that $\tilde{L}_0(n,m)=0$.

From Sakhanenko's version of the KMT inequality as stated, for example,
in Theorem~2.1 and Corollary~3.2 of \cite{Lif},
%
%e2.11 #&#
\begin{equation}
\label{E:diagramsuniform1} \mathbb{P} \biggl(Y_n^{m,1}\ge
\frac{\sqrt{n}\epsilon}{\sqrt {m}(8k-4)} \biggr)\le\bigl(1+c_2(m)\sqrt{n}\bigr) \exp
\biggl(-c_1(m)\frac{\sqrt{n}\epsilon}{\sqrt{m}(8k-4)} \biggr),
\end{equation}
where, as $m\to+\infty$, $c_1(m)\sim C_1/\sqrt{m}$ and $c_2(m)\sim
C_2/\sqrt{m}$, for
some absolute constants $C_1>0$ and $C_2>0$.
Moreover,
%
%e2.12 #&#
\begin{eqnarray}
\label{E:diagramsuniform2} \mathbb{P} \biggl(W_n^{1}\ge
\frac{\sqrt{n}\epsilon}{\sqrt{m}(8k-4)} \biggr) &\le&2n\mathbb{P} \biggl(\bigl|\tilde{B}_2^{1}\bigr|
\ge\frac{\sqrt{n}\epsilon}{\sqrt {m}(16k-8)} \biggr)
\nonumber
\\
&=&4n\mathbb{P} \biggl(\tilde{B}_2^{1}\ge
\frac{\sqrt{n}\epsilon}{\sqrt {m}(16k-8)} \biggr)
\\
&\le&4\mathrm{e} n \exp \biggl(-\frac{n\epsilon^2}{4\mathrm{e} m(16k-8)^2} \biggr).\nonumber
\end{eqnarray}
Combining \eqref{E:diagramsuniform1} and \eqref{E:diagramsuniform2},
letting $\epsilon< 1$,
and since $m(n)=\mathrm{o}(n^{1/4})$ (or simply, $m(n)=\mathrm{o}(\sqrt n)$, to get a
meaningful bound),
%
%e2.13 #&#
\begin{equation}
\label{E:diagramsmultibrown11} \mathbb{P} \bigl(\bigl\llvert \tilde{R}_{k}(n,m)-
\bigl(\tilde{L}_k(n,m)-\tilde {L}_{k-1}(n,m) \bigr) \bigr
\rrvert \ge\sqrt{mn}\epsilon \bigr)\le C_3 \sqrt{mn}\exp \biggl(-
\frac
{\sqrt{n}\epsilon}{C_{3}m} \biggr),
\end{equation}
for $1\le k \le r$, and where $C_3$ is a positive constant depending on $k$,
which for fixed $r$ can be chosen only depending on $r$.
For any $x_1\ge x_2\ge \cdots\ge x_r > 2$, $r\ge1$, and $0 < \epsilon<
\min(1,x_r - 2)$,
%
%e2.14 #&#
\begin{eqnarray}
\label{E:diagramsmultiestimate1} &&\mathbb{P} \biggl(\frac{\tilde{R}_1(n,m)}{\sqrt{mn}} \ge x_1,
\frac{\tilde{R}_2(n,m)}{\sqrt{mn}} \ge x_2,\ldots,\frac{\tilde
{R}_r(n,m)}{\sqrt{mn}}\ge
x_r \biggr)\nonumber
\\
&&\quad \le\mathbb{P} \biggl(\frac{\tilde{L}_1(n,m)-\tilde{L}_0(n,m)}{\sqrt{mn}} \ge x_1-\epsilon,\ldots,
\frac{\tilde{L}_r(n,m)-\tilde
{L}_{r-1}(n,m)}{\sqrt{mn}} \ge x_r-\epsilon \biggr)
\\
&&\qquad {} +\sum_{i=1}^{r}\mathbb{P} \biggl(
\frac{\tilde{R}_i(n,m)-(\tilde{L}_i(n,m)-\tilde
{L}_{i-1}(n,m))}{\sqrt{mn}}\ge\epsilon \biggr)
\nonumber
\end{eqnarray}
and
%
%e2.15 #&#
\begin{eqnarray}
\label{E:diagramsmultiestimate2} &&\mathbb{P} \biggl(\frac{\tilde{R}_1(n,m)}{\sqrt{mn}} \ge x_1,
\frac{\tilde{R}_2(n,m)}{\sqrt{mn}} \ge x_2,\ldots,\frac{\tilde
{R}_r(n,m)}{\sqrt{mn}}\ge
x_r \biggr)\nonumber
\\
&&\quad \ge\mathbb{P} \biggl(\frac{\tilde{L}_1(n,m)-\tilde{L}_{0}(n,m)}{\sqrt {mn}}\ge x_1+\epsilon,\ldots,
\frac{\tilde{L}_r(n,m)-\tilde
{L}_{r-1}(n,m)}{\sqrt{mn}}\ge x_r+\epsilon \biggr)
\\
&&\qquad {} -\sum_{i=1}^{r}\mathbb{P} \biggl(
\frac{(\tilde{L}_i(n,m)-\tilde{L}_{i-1}(n,m))-\tilde
{R}_i(n,m)}{\sqrt{mn}}\ge\epsilon \biggr),
\nonumber
\end{eqnarray}
with again the convention that $\tilde{L}_0(n,m)=0$.

Combining \eqref{E:brownianeigen} with Theorem~\ref{T:multi1} of the
\hyperref[app]{Appendix}, when $m$ and $n$
simultaneously converge to infinity, the large
deviations for $(\tilde{L}_k(n,m))_{1\le k \le r}$ are then given by:
%
%e2.16 #&#
%e2.17 #&#
\begin{eqnarray}
\label{E:diagramsmultibrown}&& \lim_{n \to\infty} \frac{1}{m(n)} \log\mathbb{P}
\biggl( \frac{\tilde
{L}_1(n,m(n))} {
\sqrt{m(n)n}} \ge x_1,\ldots, \frac{\tilde{L}_r(n,m(n))-\tilde{L}_{r-1}(n,m(n))} {
\sqrt{m(n)n}}\ge
x_r \biggr)\nonumber
\\[-8pt]\\[-8pt]
&&\quad =-2\sum_{i=1}^r \int_2^{x_i}
\sqrt{(z/2)^2-1}\,\mathrm{d}z,\nonumber
\end{eqnarray}
for all $x_1\ge x_2\ge\cdots\ge x_r > 2$, while for any $x<2$ and
$1\le k \le r$,
%
%e2.18 #&#
\begin{equation}
\label{cor:leftside} \lim_{n \to\infty}\frac{1}{m(n)}\log\mathbb{P}
\biggl(\frac{\tilde
{L}_k(n,m(n))-\tilde{L}_{k-1}(n,m(n))} {
\sqrt{m(n)n}}\le x \biggr)=-\infty.
\end{equation}

This implies that,
\begin{eqnarray*}
&&\mathbb{P} \biggl(\frac{\tilde{L}_1(n,m(n))}{\sqrt{m(n)n}} \ge x_1\pm\epsilon,\ldots,
\frac{\tilde{L}_r(n,m(n))-\tilde
{L}_{r-1}(n,m(n))}{\sqrt{m(n)n}} \ge x_r\pm\epsilon \biggr)
\\
&&\quad  =\exp \bigl(-m(n) \bigl(I_r(x_1\pm\epsilon,
\ldots,x_r\pm\epsilon)+\mathrm{o}(1) \bigr) \bigr),
\end{eqnarray*}
where $\mathrm{o}(1)$ indicates a quantity converging to zero as $n$ converges
to infinity.
Combining this fact with \eqref{E:diagramsmultibrown11}, for any $1\le
k \le r$,\vspace*{1pt}
%
%e2.19 #&#
\begin{eqnarray}
\label{E:diagramsmultiestimate3} &&\frac{\mathbb{P} (\llvert \tilde{R}_{k}(n,m)- (\tilde
{L}_k(n,m)-\tilde{L}_{k-1}(n,m) )\rrvert \ge\sqrt{mn}\epsilon
)}{\mathbb{P} (\tilde{L}_1(n,m) \ge\sqrt{mn}(x_1\pm\epsilon
),\ldots,\tilde{L}_r(n,m)-\tilde{L}_{r-1}(n,m)\ge\sqrt{mn}(x_r\pm
\epsilon)
)}
\nonumber
\\[1pt]
&&\quad \le C_3\sqrt{mn}\exp \biggl(-\frac{\sqrt{n}\epsilon}{C_{3}m}+m
\bigl(I_r(x_1\pm \epsilon,\ldots,x_r\pm
\epsilon)+\mathrm{o}(1)\bigr) \biggr)
\\[1pt]
&&\quad =C_3\sqrt{mn}\exp \biggl(\frac{\sqrt{n}}{m} \biggl(-
\frac{\epsilon
}{C_3}+\frac{m^2}{\sqrt{n}} \bigl(I_r(x_1\pm
\epsilon,\ldots,x_r\pm \epsilon)+\mathrm{o}(1) \bigr) \biggr) \biggr)
\longrightarrow0,\nonumber
\end{eqnarray}
as $m, n \to\infty$, $m=\mathrm{o}(n^{1/4})$. From \eqref
{E:diagramsmultiestimate1} and \eqref{E:diagramsmultiestimate3}, and
since $m=\mathrm{o}(n^{1/4})$,\vspace*{1pt}
%
%e2.20 #&#
\begin{eqnarray}\label{E:new1}
&&\limsup_{n \to\infty}\frac{1}{m(n)} \log\mathbb{P} \biggl(
\frac{\tilde
{R}_1(n,m(n))}{\sqrt{m(n)n}} \ge x_1,\ldots,\frac{\tilde
{R}_r(n,m(n))}{\sqrt{m(n)n}}\ge
x_r \biggr)\nonumber
\\[1pt]
&&\quad \le\limsup_{n\to\infty}\frac{1}{m(n)} \log2\mathbb{P} \biggl(
\frac
{\tilde{L}_1(n,m(n))}{\sqrt{m(n)n}} \ge x_1-\epsilon,\ldots,
\nonumber
\\[-8pt]\\[-8pt]
&&\hphantom{\quad \le\limsup_{n\to\infty}\frac{1}{m(n)} \log2\mathbb{P} \biggl(} \frac{\tilde{L}_r(n,m(n))-\tilde{L}_r(n,m(n))}{\sqrt{m(n)n}}\ge x_r-\epsilon \biggr)
\nonumber
\\[1pt]
&&\quad =-I_r(x_1-\epsilon,\ldots,x_r-\epsilon).
\nonumber
\end{eqnarray}
Likewise, from \eqref{E:diagramsmultiestimate2} and \eqref
{E:diagramsmultiestimate3},\vspace*{1pt}
%e2.21 #&#
\begin{eqnarray}\label{E:new2}
&&\liminf_{n \to\infty}\frac{1}{m(n)} \log\mathbb{P} \biggl(
\frac{\tilde
{R}_1(n,m(n))} {
\sqrt{m(n)n}} \ge x_1,\ldots,\frac{\tilde{R}_r(n,m(n))}{\sqrt {m(n)n}}\ge
x_r \biggr)\nonumber
\\[1pt]
&&\quad \ge\liminf_{n \to\infty}\frac{1}{m(n)} \log\frac{1}{2}
\mathbb{P} \biggl(\frac{\tilde{L}_1(n,m(n))}{\sqrt {m(n)n}} \ge x_1+\epsilon,\ldots,
\nonumber
\\[-8pt]\\[-8pt]
&&\hphantom{\quad \ge\liminf_{n \to\infty}\frac{1}{m(n)} \log\frac{1}{2}
\mathbb{P} \biggl(} \frac{\tilde{L}_r(n,m(n))-
\tilde{L}_r(n,m(n))}{\sqrt{m(n)n}}\ge x_r+\epsilon \biggr)
\nonumber
\\[1pt]
&&\quad =-I_r(x_1+\epsilon,\ldots,x_r+\epsilon).
\nonumber
\end{eqnarray}
Now letting $\epsilon\to0$,\vspace*{1pt}
\begin{eqnarray*}
&&\lim_{n \to\infty} \frac{1}{m(n)} \log\mathbb{P} \biggl(
\frac{\tilde{R}_1(n,m(n))}{\sqrt{m(n)n}} \ge x_1,\ldots,\frac{\tilde{R}_r(n,m(n))}{\sqrt{m(n)n}}\ge
x_r \biggr)
\nonumber
\\ [1pt]
&&\quad  =-2\sum_{i=1}^r \int
_2^{x_i}\sqrt{(z/2)^2-1}\,\mathrm{d}z,
\nonumber
\end{eqnarray*}
for any $x_1\ge x_2\ge \cdots\ge x_r > 2$.
Next, assume that $x_1\ge x_2\ge \cdots\ge x_k > x_{k+1}= \cdots= x_r =
2$, $1\le k\le r$,
with the convention that $k=r$ corresponds to $x_1\ge x_2 \ge\cdots\ge
x_r > 2$.
Under the conditions given in Theorem~\ref{T:diagramsmulti}, for any
$\epsilon>0$,
\begin{eqnarray*}
&&\liminf_{n \to\infty} \frac{1}{m(n)} \log\mathbb{P} \biggl(
\frac{\tilde{R}_1(n,m(n))}{\sqrt{m(n)n}} \ge x_1,\ldots,\frac{\tilde{R}_r(n,m(n))}{\sqrt{m(n)n}}\ge
x_r \biggr)
\nonumber
\\
&&\quad  \ge-2\sum_{i=1}^k \int
_2^{x_i} \sqrt{(z/2)^2-1}\,
\mathrm{d}z-2\sum_{i=k+1}^r \int
_2^{2+\epsilon}\sqrt {(z/2)^2-1}\,\mathrm{d}z.
\nonumber
\end{eqnarray*}
Letting $\epsilon\to0$, gives
%
%e2.22 #&#
\begin{eqnarray}\label{E:T:inf}
&&\liminf_{n \to\infty} \frac{1}{m(n)} \log\mathbb{P} \biggl(
\frac{\tilde{R}_1(n,m(n))}{\sqrt{m(n)n}} \ge x_1,\ldots,\frac{\tilde{R}_r(n,m(n))}{\sqrt{m(n)n}}\ge
x_r \biggr)
\nonumber
\\[-8pt]\\[-8pt]
&&\quad  \ge-2\sum_{i=1}^r \int
_2^{x_i}\sqrt{(z/2)^2-1}\,
\mathrm{d}z,\nonumber
\end{eqnarray}
while,
%
%e2.23 #&#
\begin{eqnarray}\label{E:T:sup}
&&\limsup_{n \to\infty} \frac{1}{m(n)} \log\mathbb{P} \biggl(
\frac{\tilde
{R}_1(n,m(n))}{\sqrt{m(n)n}} \ge x_1,\ldots,\frac{\tilde
{R}_r(n,m(n))}{\sqrt{m(n)n}}\ge
x_r \biggr)
\nonumber
\\
&&\quad \le\limsup_{n \to\infty} \frac{1}{m(n)} \log\mathbb{P} \biggl(
\frac{\tilde{R}_1(n,m(n))}{\sqrt{m(n)n}} \ge x_1,\ldots,\frac{\tilde{R}_k(n,m(n))}{\sqrt{m(n)n}}\ge
x_k \biggr)
\\
&&\quad = -2\sum_{i=1}^r \int
_2^{x_i}\sqrt{(z/2)^2-1}\,
\mathrm{d}z.\nonumber
\end{eqnarray}
Combining \eqref{E:T:inf} and \eqref{E:T:sup} proves \eqref
{tablemulti:finite}.\\

Fix $x < 2$ and let $0< \epsilon< \min(1,2-x)$, then
%
%e2.24 #&#
%e2.25 #&#
\begin{eqnarray}
\label{E:estimate} \mathbb{P} \biggl(\frac{\tilde{R}_k(n,m)}{\sqrt{mn}}\le x \biggr)&\le& \mathbb{P}
\biggl(\frac{\tilde{L}_k(n,m)-\tilde{L}_{k-1}(n,m)}{\sqrt {mn}}\le x+\epsilon \biggr)\nonumber
\\[-8pt]\\[-8pt]
&&{}+\mathbb{P} \biggl( \frac{|\tilde{R}_k(n,m)-(\tilde{L}_k(n,m)-\tilde
{L}_{k-1}(n,m))|}{\sqrt{mn}}\ge\epsilon \biggr),\nonumber
\end{eqnarray}
for any $1 \le k \le r$.
From \eqref{cor:leftside},
\[
\frac{1}{m}\log\mathbb{P} \biggl(\frac{\tilde{L}_k(n,m)-\tilde
{L}_{k-1}(n,m)}{\sqrt{mn}}\le x+\epsilon \biggr)
\longrightarrow-\infty,
\]
and from \eqref{E:diagramsmultibrown11},
%
%e2.26 #&#
\begin{eqnarray}
&&\frac{1}{m}\log\mathbb{P} \biggl(\frac{|\tilde{R}_{k}(n,m)-(\tilde
{L}_k(n,m)-\tilde{L}_{k-1}(n,m))|}{\sqrt{mn}} \ge\epsilon \biggr)
\nonumber
\\[-8pt]\\[-8pt]
&&\quad  \le\frac{\log(C_3\sqrt{mn})}{m}-\frac{\sqrt{n} \epsilon}{C_3 m^2} \longrightarrow-\infty,\nonumber
\end{eqnarray}
as $m, n \to\infty$, $m=\mathrm{o}(n^{1/4})$. Thus for any $x<2$ and $1\le k
\le r$,
%
%e2.27 #&#
\begin{equation}
\lim_{n \to\infty}\frac{1}{m(n)}\log\mathbb{P} \biggl(
\frac{\tilde{R}_k(n,m(n))} {
\sqrt{m(n)n}}\le x \biggr)=-\infty,
\end{equation}
which proves \eqref{tablemulti:infinity} in Theorem~\ref{T:diagramsmulti}.
\begin{pf*}{Proof of Theorem~\ref{P:tableautraceless}}
First, \eqref{th2_zero} is just a direct consequence of \eqref
{tablemulti:finite}.
Next, we prove \eqref{th2_K}. Fix $y_1\ge y_2 \ge\cdots\ge y_r $,
with $y_r<2$.
If $K(y_r)<+\infty$, then there exists $\delta>0$ such that
$K(y_r-\delta)<+\infty$ and such that for any $0< \epsilon< \min(1,
\delta,2-y_r)$,
%
%e2.28 #&#
\begin{eqnarray}
\label{E:diagramsmultiestimate4} &&\mathbb{P} \biggl(\frac{\tilde{R}_1(n,m)}{\sqrt{mn}} \le y_1,\ldots,
\frac
{\tilde{R}_r(n,m)}{\sqrt{mn}}\le y_r \biggr)
\nonumber
\\
&&\quad \le\mathbb{P} \biggl(\frac{\tilde{L}_1(n,m)-\tilde{L}_0(n,m)}{\sqrt {mn}} \le y_1+\epsilon,\ldots,
\frac{\tilde{L}_r(n,m)-\tilde
{L}_{r-1}(n,m)}{\sqrt{mn}}\le y_r+\epsilon \biggr)
\\
&&\qquad {} +\sum_{i=1}^{r}\mathbb{P} \biggl(
\frac{|\tilde{R}_i(n,m)-(\tilde
{L}_i(n,m)-\tilde{L}_{i-1}(n,m))|}{\sqrt{mn}}\ge\epsilon \biggr)
\nonumber
\end{eqnarray}
and
%
%e2.29 #&#
\begin{eqnarray}
\label{E:diagramsmultiestimate5} &&\mathbb{P} \biggl(\frac{\tilde{R}_1(n,m)}{\sqrt{mn}} \le y_1,\ldots,
\frac
{\tilde{R}_r(n,m)}{\sqrt{mn}}\le y_r \biggr)
\nonumber
\\
&&\quad \ge\mathbb{P} \biggl(\frac{\tilde{L}_1(n,m)-\tilde{L}_0(n,m)}{\sqrt {mn}} \le y_1-\epsilon,\ldots,
\frac{\tilde{L}_r(n,m)-\tilde
{L}_{r-1}(n,m)}{\sqrt{mn}}\le y_r-\epsilon \biggr)
\\
&&\qquad {} -\sum_{i=1}^{r}\mathbb{P} \biggl(
\frac{|\tilde{R}_i(n,m)-(\tilde
{L}_i(n,m)-\tilde{L}_{i-1}(n,m))|}{\sqrt{mn}}\ge\epsilon \biggr),
\nonumber
\end{eqnarray}
with once more the convention that $\tilde{L}_0(n,m)=0$.

Combining \eqref{E:brownianeigen} with Corollary~\ref{coro2}, when
$m$ and $n$ simultaneously
converge to infinity,
%
%e2.30 #&#
%e2.31 #&#
\begin{eqnarray}
\label{E:diagramsmultibrown111} &&\lim_{n \to\infty} \frac{1}{m(n)^2} \log\mathbb{P}
\biggl( \frac{\tilde{L}_1(n,m(n))}{\sqrt{m(n)n}} \leq y_1,\ldots, \frac{\tilde{L}_r(n,m(n))-\tilde{L}_{r-1}(n,m(n))} {
\sqrt{m(n)n}} \leq
y_r \biggr)\nonumber
\\ [-8pt]\\[-8pt]
&&\quad =-K(y_r),\nonumber
\end{eqnarray}
for all $y_r\le y_{r-1} \le\cdots\le y_1$ with $y_r<2$. Thus
\begin{eqnarray*}
&&\mathbb{P} \biggl(\frac{\tilde{L}_1(n,m(n))}{\sqrt{m(n)n}} \le y_1\pm \epsilon,\ldots,
\frac{\tilde{L}_r(n,m(n))-\tilde{L}_{r-1}(n,m(n))}{\sqrt {m(n)n}}\le y_r\pm\epsilon \biggr)
\nonumber
\\
&&\quad =\exp \bigl(-m(n)^2 \bigl(K(y_r\pm\epsilon)+
\mathrm{o}(1) \bigr) \bigr),
\nonumber
\end{eqnarray*}
where $\mathrm{o}(1)$ is meant for an expression converging to zero as
$n$ converges to infinity.
Combining this last fact with \eqref{E:diagramsmultibrown11}, for any
$1\le k \le r$,
%
%e2.32 #&#
%e2.33 #&#
\begin{eqnarray}
\label{E:diagramsmultiestimate6} &&\frac{\mathbb{P} (\llvert \tilde{R}_{k}(n,m)- (\tilde
{L}_k(n,m)-\tilde{L}_{k-1}(n,m) )\rrvert \ge\sqrt{mn}\epsilon
)}{\mathbb{P} (\tilde{L}_1(n,m)
\le\sqrt{mn}(y_1\pm\epsilon),\ldots,\tilde{L}_r(n,m)-\tilde
{L}_{r-1}(n,m)\le\sqrt{mn}(y_r\pm\epsilon)
)}\nonumber
\\[-8pt]\\[-8pt]
&&\quad \le C_3\sqrt{mn}\exp \biggl\{\frac{\sqrt{n}}{m} \biggl(-
\frac{\epsilon
}{C_3}+\frac{m^3}{\sqrt{n}} \bigl(K(y_r\pm\epsilon)+
\mathrm{o}(1) \bigr) \biggr) \biggr\}\longrightarrow0,\nonumber
\end{eqnarray}
as $m, n \to\infty$, $m=\mathrm{o}(n^{1/6})$. Repeating previous arguments, letting
$\epsilon\to0$, and since $m=\mathrm{o}(n^{1/6})$,
%
%e2.34 #&#
\begin{equation}
\label{E:propo1} \lim_{n \to\infty} \frac{1}{m(n)^2} \log\mathbb{P}
\biggl(\frac{\tilde{R}_1(n,m(n))}{\sqrt{m(n)n}} \le y_1,\ldots,\frac{\tilde{R}_r(n,m(n))}{\sqrt{m(n)n}}\le
y_r \biggr)=-K(y_r),
\end{equation}
for $y_r\le y_{r-1} \le\cdots\le y_1 $, with $y_r<2$ and
$K(y_r)<+\infty$.

Now for fixed $y_1\ge y_2 \ge\cdots\ge y_r $, $y_r<2$, let us tackle
the case $K(y_r)=+\infty$.
Then,
%
%e2.35 #&#
\begin{eqnarray}
\label{T1:1} &&\mathbb{P} \biggl(\frac{\tilde{R}_1(n,m)}{\sqrt{mn}} \le y_1,\ldots,
\frac
{\tilde{R}_r(n,m)}{\sqrt{mn}}\le y_r \biggr)
\nonumber
\\
&&\quad  \le\mathbb{P} \biggl(\frac{\tilde{L}_r(n,m)-\tilde{L}_{r-1}(n,m)}{\sqrt {mn}}\le y_r+\epsilon \biggr)
\\
&&\qquad {} + \mathbb{P} \biggl(\frac{|\tilde{R}_r(n,m)-(\tilde{L}_r(n,m)-\tilde
{L}_{r-1}(n,m))|}{\sqrt{mn}}\ge\epsilon \biggr).\nonumber
\end{eqnarray}
As $m$ and $n$ simultaneously converge to infinity with $m=\mathrm{o}(n^{1/6})$,
the second term on the right of \eqref{T1:1} is exponentially negligible
with speed $m^2$, that is,
%
%e2.36 #&#
\begin{eqnarray}
&&\frac{1}{m^2}\log\mathbb{P} \biggl(\frac{|\tilde{R}_{k}(n,m)-(\tilde
{L}_k(n,m)-\tilde{L}_{k-1}(n,m))|}{\sqrt{mn}} \ge\epsilon \biggr)
\nonumber
\\[-8pt]\\[-8pt]
&&\quad  \le\frac{\log(C_3\sqrt{mn})}{m^2}-\frac{\sqrt{n} \epsilon}{C_3 m^3} \longrightarrow-\infty,\nonumber
\end{eqnarray}
while the first term is,
from \eqref{E:diagramsmultibrown111}, dominated by
$\mathrm{e}^{-m(n)^2K(y_r+\epsilon)}$.
Thus \eqref{E:propo1}, in this case, follows by letting $\epsilon\to0$.

Now let $2=y_r\le y_{r-1} \le\cdots\le y_1$, then for any $\epsilon>0$,
%
%e2.37 #&#
\begin{eqnarray}
&&\liminf_{n \to\infty} \frac{1}{m(n)^2} \log\mathbb{P} \biggl(
\frac{\tilde{R}_1(n,m(n))}{\sqrt{m(n)n}} \le y_1,\ldots,\frac{\tilde{R}_r(n,m(n))}{\sqrt{m(n)n}}\le
y_r \biggr)
\nonumber
\\
&&\quad \ge\liminf_{n \to\infty} \frac{1}{m(n)^2} \log\mathbb{P} \biggl(
\frac
{\tilde{R}_1(n,m(n))}{\sqrt{m(n)n}} \le y_1,\ldots,\frac{\tilde
{R}_r(n,m(n))}{\sqrt{m(n)n}}\le2-\epsilon
\biggr)
\\
&&\quad =-K(2-\epsilon).
\nonumber
\end{eqnarray}
Again, letting $\epsilon\to0$, and since $K$ is continuous (see the
\hyperref[app]{Appendix} for a proof),
\[
\liminf_{n \to\infty} \frac{1}{m(n)^2} \log\mathbb{P} \biggl(
\frac
{\tilde{R}_1(n,m(n))}{\sqrt{m(n)n}} \le y_1,\ldots,\frac{\tilde
{R}_r(n,m(n))}{\sqrt{m(n)n}}\le
y_r \biggr)\ge-K(2)=0.
\]
Clearly,
\[
\limsup_{n \to\infty} \frac{1}{m(n)^2} \log\mathbb{P} \biggl(
\frac{\tilde{R}_1(n,m(n))}{\sqrt{m(n)n}} \le y_1,\ldots,\frac{\tilde{R}_r(n,m(n))}{\sqrt{m(n)n}}\le
y_r \biggr)\le0,
\]
which proves the case $y_r=2$, and finishes the proof of the first part of
Theorem~\ref{P:tableautraceless}. Lemma~\ref{lemma_appendix} of the \hyperref[app]{Appendix}
gives a proof of \eqref{ratefunction_traceless}.

When $x\le0$, $\mathcal{M}_{0}((-\infty,x])$ is empty so $K(x)=+\infty
$ and
when $x\ge2$, the semicircular probability measure belongs to $\mathcal
{M}_{0}((-\infty,x])$,
thus $K(x)=0$. When $0<x<2$, the closed form expression of $K$ given via
\eqref{ratefunction_explicit} can indeed be derived using the
techniques developed
in \cite{NMV}.
Denote by $\mu_0$ the zero mean probability measure supported on
$(-\infty,x]$,
minimizing
%
%e2.38 #&#
\begin{equation}
I(\mu)=\frac{1}{2}\int y^2 \mu(\mathrm{d}y)-\iint\log|t-y|
\mu(\mathrm{d}t)\mu(\mathrm{d}y)-\frac{3}{4}
\end{equation}
(the existence and uniqueness of $\mu_0$ follows from Theorem~1.3 of
Chapter~1 of \cite{SF}.
Moreover, in view of Theorem~2.5 of Chapter IV of \cite{SF},
$\mu_0$ is absolutely continuous with continuous density
$\rho_0$, while from Theorem~1.10 and Theorem~1.11 of Chapter IV in
\cite{SF},
its support is a finite interval).
Let us now proceed to explicitly find $\rho_0$. To do so, consider the
Lagrange function
\[
E(\mu)=I(\mu)+c_1 \biggl(\int\mu(\mathrm{d}y)-1
\biggr)+c_2\int y \mu(\mathrm{d}y),
\]
where the Lagrange multipliers $c_1$ and $c_2$ correspond to the
constraints that $\mu$ is a
zero mean probability measure. Let $[L',x]$ be the support of $\rho_0$,
and for any continuous function $h$ supported on $[L',x]$ such that
$h(y)\ge-\rho_0(y)$, let
%
%e2.39 #&#
%e2.40 #&#
%e2.41 #&#
\begin{eqnarray}
E (\rho_0+\epsilon h)&=&\frac{1}{2}\int y^2 \bigl(
\rho_0(y)+\epsilon h(y) \bigr)\,\mathrm{d}y\nonumber
\\
&&{}-\iint\log|t-y| \bigl(\rho_0(t)+\epsilon h(t) \bigr) \bigl(\rho
_0(y)+\epsilon h(y) \bigr)\,\mathrm{d}t \,\mathrm{d}y-
\frac{3}{4}
\\
&&{}+c_1 \biggl(\int \bigl(\rho_0(y)+\epsilon h(y) \bigr)\,
\mathrm{d}y-1 \biggr)+c_2\int y \bigl(\rho_0(y)+\epsilon
h(y) \bigr) \,\mathrm{d}y.\nonumber
\end{eqnarray}
Thus
\[
\frac{\mathrm{d} E (\rho_0+\epsilon h)}{\mathrm{d} \epsilon} \biggl|_{\epsilon=0}=0,
\]
gives
%
%e2.42 #&#
\begin{equation}
\label{equality} \int \biggl(\frac{y^2}{2}-2\int\log|t-y|\rho_0(t)
\,\mathrm{d}t+c_1+c_2 y \biggr)h(y)\,\mathrm{d}y=0,
\end{equation}
for any continuous function $h$ such that $h(y)\ge-\rho_0(y)$. Let
\[
g(y)=\frac{y^2}{2}-2\int\log|t-y|\rho_0(t)\,
\mathrm{d}t+c_1+c_2 y,
\]
which is a continuous function on $[L',x]$. Let $h(y)=g^+(y)$, then
\eqref{equality} yields
\[
\int_{g(y)\ge0} g(y)^2 \,\mathrm{d}y =0,
\]
thus $g(y)\le0$ for $y\in[L',x]$. Likewise, letting
%
%e2.43 #&#
\begin{equation}
h(y)= %
\cases{ 0, &\quad \mbox{if $g(y)> 0$,}
\cr
g(y), &\quad \mbox{if $-
\rho_0(y)\le g(y)\le0$,}
\cr
-\rho_0(y), &\quad \mbox{if $g(y)<
-\rho_0(y)$,} } %
\end{equation}
then \eqref{equality} yields
$g(y)\ge0$ for $y\in[L',x]$. Thus,
%
%e2.44 #&#
\begin{equation}
\label{prop_equ2} \frac{y^2}{2}-2\int\log|t-y|\rho_0(t)\,
\mathrm{d}t+c_1+c_2 y=0,
\end{equation}
for any $y \in[L',x]$.
In turn, differentiating \eqref{prop_equ2} with respect to $y$ further gives,
%
%e2.45 #&#
\begin{equation}
y-2 \mbox{ p.v. }\int\frac{\rho_0(t)}{y-t}\,\mathrm{d}t+c_2=0,
\end{equation}
where $\mbox{ p.v. }$ denotes the Cauchy principal value.

Let $L=L'-x$ and $f_x(t)=\rho_0(t+x)$ be supported on $[L,0]$,
then the finite Hilbert transform
\[
\frac{1}{\uppi}\mbox{ p.v. }\int_{L'}^{x}
\frac{\rho_0(t)}{y-t}\,\mathrm{d}t=\frac
{y+c_2}{2\uppi},
\]
becomes
\[
\frac{1}{\uppi}\mbox{ p.v. }\int_{L}^{0}
\frac{f_x(t)}{y-t}\,\mathrm{d}t=\frac
{x+y+c_2}{2\uppi},
\]
for any $y \in[L,0]$. From Section~4.3 of \cite{Tri}, this finite
Hilbert transform
can be inverted as
%
%e2.46 #&#
\begin{equation}
\label{prop_equ1} f_x(y)=\frac{1}{\uppi\sqrt{(y-L)(-y)}} \biggl(\mbox{p.v. }\int
_L^0 \frac{\sqrt{(t-L)(-t)}}{t-y} \frac{x+t+c_2}{2\uppi}\,
\mathrm{d}t+c_3 \biggr),
\end{equation}
where $L\le y \le0$. Moreover,
%
%e2.47 #&#
%e2.48 #&#
\begin{eqnarray}
&&\mbox{p.v. }\int_L^0 \frac{\sqrt{(t-L)(-t)}}{t-y}
\frac{x+t+c_2}{2\uppi
}\,\mathrm{d}t\nonumber
\\[-8pt]\\[-8pt]
&&\quad =\frac{1}{16} \bigl(4 c_2 (L-2 y)+L^2+4 L (x+y)-8
y (x+y) \bigr).\nonumber
\end{eqnarray}
Since $f_x(L)=0$,
\[
c_3=\tfrac{1}{16} \bigl(4L(c_2+x)+3
L^2 \bigr),
\]
which when plugged into \eqref{prop_equ1} yields
\[
f_x(y)=\frac{\sqrt{y (L-y)} (2 c_2+L+2 (x+y))}{4 \uppi y}.
\]
Now, the two constraints $\int \mathrm{d}\mu_0(y)=1$ and $\int y\,\mathrm{d}\mu_0(y)=0$, yield
\[
\int_L^0 y f_x(y) \,\mathrm{d}y
+x=0, \qquad \int_L^0 f_x(y)\,
\mathrm{d}y=1,
\]
leading to
%
%e2.49 #&#
\begin{equation}
L=\frac{2\, 2^{2/3}  (\sqrt{81 x^2+12}-9 x )^{2/3}-4\, 6^{1/3}} {
3^{2/3}  (\sqrt{81 x^2+12}-9 x )^{1/3}}
\end{equation}
and
%
%e2.50 #&#
%e2.51 #&#
\begin{eqnarray}
c_2&=&\frac{2\, 3^{2/3}-\sqrt[3]{6}
(\sqrt{81 x^2+12}-9 x )^{2/3}}{2^{2/3}  (\sqrt{81
x^2+12}-9 x )^{1/3}}\nonumber
\\[-8pt]\\[-8pt]
&&{}-\frac{\sqrt[3]{2}\
3^{2/3}  (\sqrt{81 x^2+12}-9 x )^{2/3}
+\sklfrac{6\, 2^{2/3} \sqrt[3]{3}}{ (\sqrt{81 x^2+12}-9 x )^{2/3}}+6}{18
x}-x.\nonumber
\end{eqnarray}
Integrate \eqref{prop_equ2} with respect to $\mu_0$ to get
\[
\iint\log|y-t|\mu_0(\mathrm{d}t)\mu_0(\mathrm{d}y)=
\frac{1}{4}\int y^2 \mu_0(\mathrm{d}y)+
\frac{c_1}{2},
\]
while $c_1$ is determined by substituting $y=x$ in \eqref{prop_equ2},
\[
c_1=-\frac{x^2}{2}+2\int\log|x-t|\mu_0(
\mathrm{d}t)-c_2 x.
\]
Finally,
%
%e2.52 #&#
\begin{eqnarray}\label{prop_new_equ1}
I(\mu_0)&=&\frac{1}{2}\int y^2 \mu_0(
\mathrm{d}y)-\iint\log|t-y|\mu_0(\mathrm{d}t)\mu _0(
\mathrm{d}y)-\frac{3}{4}
\nonumber
\\[-8pt]\\[-8pt]
&=&\frac{1}{4}\int_L^0(x+y)^2
f_x(y) \,\mathrm{d}y-\int_L^0
\log(-y) f_x(y)\,\mathrm{d}y +\frac
{x^2}{4}+
\frac{c_2x}{2}-\frac{3}{4}.\nonumber
\end{eqnarray}
Finally, inserting the above values of $L$ and $c_2$ into \eqref{prop_new_equ1}
provides the closed form expression for $K$.
\end{pf*}

%s3 #&#
\section{Proof of Theorem \texorpdfstring{\protect\ref{T:nonuniform}}{1.3} and 
Theorem \texorpdfstring{\protect\ref{P:nonuniform}}{1.4}}
Recall, see \eqref{def}, that
\[
R_1(n,m)=V_{1}(n,m)=\sup_{0=l_{0} \le l_{1} \le\cdots \le l_{m}=n} \sum
_{j=1}^{m}\bigl(S_{l_j}^{m,j}-S_{l_{j-1}}^{m,j}
\bigr)
\]
and let
\[
V_{1}^{\prime}(n,m)=\mathop{\sup_{0=l_{0} \le l_{1} \le\cdots \le
l_{m}=n}}_{
l_{j-1}=l_j \mathrm{\ for\ } j \notin J(m)}
\sum_{j=1}^{m}\bigl(S_{l_j}^{m,j}-S_{l_{j-1}}^{m,j}
\bigr),
\]
where $J(m)=\{j\dvtx p_j^m=p_{\mathrm{max}}^m\}$ as defined before. Then, Lemma~9 in
\cite{BH} asserts that:
%
%e3.1 #&#
\begin{equation}
\label{E:nonuniform221} \mathbb{E}\bigl\llvert V_1(n,m)-V_1^{\prime}(n,m)
\bigr\rrvert \le Cnp_{\mathrm{2nd}}^m,
\end{equation}
where $p_{\mathrm{2nd}}$ is the second highest probability and $C>0$ some
absolute constant.

In order to prove Theorem~\ref{T:nonuniform}, we first need a lemma.
%
%le3.1 #&#
\begin{lemma}\label{lemma2}
Let $k(m(n))$ converge to infinity with $n$ in such a way that
$k(m(n))^3/p_{\mathrm{max}}^m=\mathrm{o}(n)$, then for any $x\ge2$,
%
%e3.2 #&#
\begin{equation}
\label{E:nonuniformlemma} \lim_{n \to\infty}\frac{1}{k(m(n))} \log\mathbb{P}
\biggl(\frac
{V_1^{\prime}(n,m(n))-np_{\mathrm{max}}^m}{\sqrt{nk(m(n))p_{\mathrm{max}}^m}} \ge x \biggr)=-2 \int_2^{x}
\sqrt{(z/2)^2-1}\,\mathrm{d}z,
\end{equation}
and for any $x<2$,
%
%e3.3 #&#
\begin{equation}
\label{E:nonuniformlemma1} \lim_{n \to\infty}\frac{1}{k(m(n))} \log\mathbb{P}
\biggl(\frac
{V_1^{\prime}(n,m(n))-np_{\mathrm{max}}^m}{\sqrt{nk(m(n))p_{\mathrm{max}}^m}} \le x \biggr)=-\infty.
\end{equation}
\end{lemma}
\begin{pf} As in the proof of Theorem~\ref{T:diagramsmulti},
for any $j\in J(m)$, set $\tilde{X}_{i,j}^m=(X_{i,j}^m-p_{\mathrm{max}}^m)/\sigma_m$,
where $\sigma_m^2=p_{\mathrm{max}}^m(1-p_{\mathrm{max}}^m)$, and
set $\tilde{S}_{\ell}^{m,j}=\sum_{i=1}^{\ell}\tilde{X}_{i,j}^m$. Hence,
\[
\frac{V_1^{\prime}(n,m)-np_{\mathrm{max}}^m}{\sqrt{nk(m(n))p_{\mathrm{max}}^m}}=\Bigl(\sqrt {1-p_{\mathrm{max}}^m}
\Bigr) \frac{ \tilde{V}_1^{\prime}(n,m)}{\sqrt{nk(m(n))}},
\]
with the obvious notation for $\tilde{V}_1^{\prime}(n,m)$.
Since $k(m(n))p_{\mathrm{max}}^m\le1$, as $n \to\infty$, $p_{\mathrm{max}}^m \to0$, so
\eqref{E:nonuniformlemma} can be reduced to,
%
%e3.4 #&#
\begin{equation}
\label{E:nonuniform24} \lim_{n \to\infty}\frac{1}{k(m(n))} \log\mathbb{P}
\biggl(\frac{ \tilde{V}_1^{\prime}(n,m(n))}{\sqrt{nk(m(n))}} \ge x \biggr)=-I_1(x),
\end{equation}
for any $x\ge2$. Moreover, \eqref{E:nonuniformlemma1} can be reduced to,
%
%e3.5 #&#
\begin{equation}
\label{E:nonuniform25} \lim_{n \to\infty}\frac{1}{k(m(n))} \log\mathbb{P}
\biggl(\frac{ \tilde{V}_1^{\prime}(n,m(n))}{\sqrt {nk(m(n))}} \le x \biggr)=-\infty,
\end{equation}
for any $x<2$. Now,
%
%e3.6 #&#
\begin{equation}
\tilde{V}_{1}^{\prime}(n,m)=\mathop{\sup_{0=l_{0} \le l_{1} \le
\cdots \le l_{m}=n}}_{
l_{j-1}=l_j \mathrm{\ for\ } j \notin J(m)}
\sum_{j=1}^{m}\bigl(\tilde
{S}_{l_j}^{m,j}-\tilde{S}_{l_{j-1}}^{m,j}
\bigr),
\end{equation}
where
%
%e3.7 #&#
\begin{equation}
\Cov\bigl(\tilde{S}_{\ell}^{m,i},\tilde{S}_{\ell}^{m,j}
\bigr)= %
\cases{ \ell, &\quad \mbox{if $i=j$,}
\cr
\rho_1\ell, &\quad
\mbox{otherwise,} } %
\end{equation}
where $\rho_1=-p_{\mathrm{max}}^m/(1-p_{\mathrm{max}}^m)$.
From its very definition, $\tilde{V}_{1}^{\prime}(n,m)$ only depends on
$(\tilde{S}_{\ell}^{m,j})_{j \in J(m)}$
and can thus be approximated, via KMT, by the Brownian functional
$F(n,k)$ with $k=\card (J(m))$ (from here onward, $m$ is short for
$m(n)$ and $k$ is short for
$k(m(n))$), where
%
%e3.8 #&#
\begin{equation}
\label{E:newbrownian} F(n,k)=\sup_{0=t_{0} \le t_{1} \le\cdots
\le t_{k}=n} \sum
_{r=1}^{k}\bigl(\tilde{B}_{t_r}^{(r)}-
\tilde{B}_{t_{r-1}}^{(r)}\bigr),
\end{equation}
where $(\tilde{B}^{(r)})_{1 \le r \le k}$ is a centered $k$-dimensional
Brownian motion
with covariance matrix
\[
\pmatrix{ 1& \rho_1 & \cdots& \rho_1
\cr
\rho_1 & 1 & \cdots& \rho_1
\cr
\vdots&\vdots&\ddots&
\vdots
\cr
\rho_1 & \rho_1 &\cdots&1 } %
t.
\]
Moreover,
%
%e3.9 #&#
\begin{equation}
\label{E:neweq} F(n,k)\stackrel{\mathcal{L}} {=}\sqrt{n}F(1,k),
\end{equation}
while from Corollary~3.2 and Corollary~3.3 in \cite{HL1},
%
%e3.10 #&#
\begin{eqnarray}
\label{E:nonuniform22} \sqrt{1-p_{\mathrm{max}}^m} F(1,k)&
\stackrel{\mathcal{L}} {=}& \frac{\sqrt{1-kp_{\mathrm{max}}^m}-1}{k}\sum_{j=1}^k
B^{(j)}_1\nonumber
\\[-8pt]\\[-8pt]
&&{} + \sup_{0=t_{0} \le t_{1} \le\cdots \le t_{k}=1} \sum_{r=1}^{k}
\bigl(B_{t_r}^{(r)}-B_{t_{r-1}}^{(r)}\bigr) ,
\nonumber
\end{eqnarray}
where $(B^{(j)})_{1\le j \le k}$ is a standard $k$-dimensional Brownian
motion on $[0,1]$.
The first weighted sum in \eqref{E:nonuniform22} is a
Gaussian random variable with variance at most $1/k$, and
as well known (see the introductory section and the cited references therein):
%
%e3.11 #&#
\begin{equation}
\label{E:eigenbrownian} \sup_{0=t_{0} \le t_{1} \le\cdots \le t_{k}=1} \sum_{r=1}^{k}
\bigl(B_{t_r}^{(r)}-B_{t_{r-1}}^{(r)}\bigr)
\stackrel{\mathcal {L}} {=}\lambda_1^k ,
\end{equation}
where $\lambda_1^k$ is the largest eigenvalue of
a $k\times k$ element of the GUE. Next,
since $\lambda_1^k/\sqrt{k}$ satisfies a LDP with rate function $I_1$
and since the
contribution of the Gaussian term is negligible, we get via Theorem~\ref
{T:multi1} of the \hyperref[app]{Appendix}:
%
%e3.12 #&#
\begin{equation}
\label{E:nonuniform23} \lim_{k \to\infty}\frac{1}{k} \log\mathbb{P}
\bigl(F(1,k) \ge\sqrt{k}x\bigr)=-I_1(x).
\end{equation}
Now, as in the proof of Theorem~\ref{T:diagramsmulti},
%
%e3.13 #&#
%e3.14 #&#
\begin{eqnarray}\label
{E:nonuniform-concentrationback}
&&\mathbb{P} \bigl(\bigl\llvert \tilde{V}_{1}^{\prime}(n,m)-F(n,k)
\bigr\rrvert \ge\sqrt {nk}\epsilon \bigr)\nonumber\\[-8pt]\\[-8pt]
&&\quad \le k\mathbb{P} \biggl(Y_n^{m,l}\ge \frac{\sqrt{n}\epsilon}{4\sqrt{k}}
\biggr)+k\mathbb{P} \biggl(W_n^{l} \ge\frac{\sqrt{n}\epsilon}{4\sqrt{k}}
\biggr),\nonumber
\end{eqnarray}
where $l$ is any element of $J(m)$ and where
\[
Y_n^{m,l}=\max_{1\le i \le n}\bigl|
\tilde{S}_{i}^{m,l}- \tilde{B}_{i}^{(l)}\bigr|
\quad \mbox{and}\quad  W_n^l=\mathop{\sup_{0\le s,t \le
n}}_{ |s-t|\le1}
\bigl|\tilde{B}_{s}^{(l)}-\tilde{B}_{t}^{(l)}\bigr|.
\]
As in getting \eqref{E:diagramsuniform1},
%
%e3.15 #&#
\begin{equation}
\label{E:diagramsnonuniform1} \mathbb{P} \biggl(Y_n^{m,1}\ge
\frac{\sqrt{n}\epsilon} {
4\sqrt{k}} \biggr)\le\bigl(1+c_2\bigl(p_{\mathrm{max}}^m
\bigr)\sqrt{n}\bigr) \exp \biggl(-c_1\bigl(p_{\mathrm{max}}^m
\bigr)\frac{\sqrt{n}\epsilon}{4\sqrt{k}} \biggr),
\end{equation}
where $c_1(p_{\mathrm{max}}^m)\sim C_1\sqrt{p_{\mathrm{max}}^m}$ and
$c_2(p_{\mathrm{max}}^m)\sim C_2\sqrt{p_{\mathrm{max}}^m}$, for some constants $C_1$ and $C_2$,
while from \eqref{E:diagramsuniform2},
%
%e3.16 #&#
\begin{equation}
\label{E:diagramsnonuniform2} \mathbb{P} \biggl(W_n^{1}\ge
\frac{\sqrt{n}\epsilon} {
4 \sqrt{k}} \biggr) \le C_3 n \exp \biggl(-\frac{n\epsilon^2}{C_3 k}
\biggr),
\end{equation}
for some positive constant $C_3$.
Combining \eqref{E:diagramsnonuniform1} and \eqref
{E:diagramsnonuniform2}, letting $\epsilon< 1$, and since
$k(m(n))^3/\allowbreak p_{\mathrm{max}}^m=\mathrm{o}(n)$ (or simply, $k(m(n))/p_{\mathrm{max}}^m=\mathrm{o}(n)$, to get
a meaningful bound),
%
%e3.17 #&#
\begin{equation}
\label{E:diagramsnonuniformtotal} \mathbb{P} \bigl(\bigl\llvert \tilde{V}_{1}^{\prime}(n,m)-F(n,k)
\bigr\rrvert \ge\sqrt {nk}\epsilon \bigr) \le C_4k\sqrt
{np_{\mathrm{max}}^m}\exp \biggl(-\frac{\sqrt{np_{\mathrm{max}}^m}\epsilon
}{C_4\sqrt{k}}
\biggr),
\end{equation}
for some positive constant $C_4$.
From \eqref{E:nonuniform23}, for any $x>2$ and $0 < \epsilon< \min(1,x-2)$,
%
%e3.18 #&#
\begin{equation}
\label{E:nonuniformexpo} \mathbb{P}\bigl(F(n,k) \ge\sqrt{nk}(x\pm\epsilon)\bigr)=\exp\bigl
\{-k\bigl(I_1(x\pm \epsilon)+\mathrm{o}(1)\bigr)\bigr\}.
\end{equation}
Hence,
\begin{eqnarray*}
&&\frac{\mathbb{P} (\llvert \tilde{V}_{1}^{\prime}(n,m)-F(n,k)\rrvert \ge\sqrt{nk}\epsilon )}{\mathbb{P}(F(n,k) \ge\sqrt{nk}(x\pm
\epsilon))}
\nonumber
\\
&&\quad \le C_4 k\sqrt{np_{\mathrm{max}}^m}\exp
\biggl(\sqrt{\frac{np_{\mathrm{max}}^m}{k}} \biggl(-\frac{\epsilon}{C_4}+\sqrt
{\frac{k^3}{np_{\mathrm{max}}^m}} \bigl(I_1(x\pm \epsilon)+\mathrm{o}(1)
\bigr) \biggr) \biggr) \longrightarrow0,
\nonumber
\end{eqnarray*}
since $k(m(n))^3/p_{\mathrm{max}}^m=\mathrm{o}(n)$, and, again, as in the proof of
Theorem~\ref{T:diagramsmulti}, this leads to
\eqref{E:nonuniform24} for any $x>2$. Arguments similar to those
developed at the end of the proof of Theorem~\ref{T:diagramsmulti}
show that \eqref{E:nonuniform24} is valid for any $x\ge2$.

The proof of \eqref{E:nonuniform25} is similar to the uniform case.
First, from \eqref{E:nonuniform22} and
\eqref{E:eigenbrownian}, for any fixed $x<2$,
%
%e3.19 #&#
\begin{equation}
\label{E:nonuniform26} \lim_{k \to\infty}\frac{1}{k} \log\mathbb{P}
\bigl(F(1,k) \le\sqrt {k}x\bigr)=-\infty.
\end{equation}
Moreover, for any $0<\epsilon<\min(1,2-x)$,
%
%e3.20 #&#
%e3.21 #&#
\begin{eqnarray}
&&\mathbb{P} \bigl(\tilde{V}_{1}^{\prime}(n,m) \le\sqrt{nk}x
\bigr)\nonumber
\\[-8pt]\\[-8pt]
&&\quad \le\mathbb{P} \bigl(F(n,k)\le\sqrt{nk}(x+\epsilon) \bigr)+ \mathbb{P} \bigl(\bigl
\llvert \tilde{V}_{1}^{\prime}(n,m)-F(n,k)\bigr\rrvert \ge\sqrt
{nk}\epsilon \bigr),\nonumber
\end{eqnarray}
while $\mathbb{P} (\llvert \tilde{V}_{1}^{\prime}(n,m)-F(n,k)\rrvert \ge\sqrt{nk}\epsilon )$ is exponentially negligible with speed
$k(m)$. Therefore, \eqref{E:nonuniform25} holds true under the condition
$k(m(n))^3/p_{\mathrm{max}}^m=\mathrm{o}(n)$.
\end{pf}
\begin{pf*}{Proof of Theorem~\ref{T:nonuniform}}
First, so as not to further burden the notations, below $m$ is short
for $m(n)$ and $k$ is short for $k(m(n))$.
Next, set $X=(V_1(n,m)-np_{\mathrm{max}}^m)/\sqrt{nkp_{\mathrm{max}}^m}$,
$Y=(V_1(n,m)-V_1^{\prime}(n,m))/\sqrt{nkp_{\mathrm{max}}^m}$
and $Z=(V_1^{\prime}(n,m)-np_{\mathrm{max}}^m)/\sqrt{nkp_{\mathrm{max}}^m}$.
Then, for any $x>2$ and $0<\epsilon<x-2$,
%
%e3.22 #&#
\begin{equation}
\mathbb{P}(X \ge x)\le\mathbb{P}(Z \ge x-\epsilon)+\mathbb{P}\bigl(|Y| \ge \epsilon\bigr)
\end{equation}
and
%
%e3.23 #&#
\begin{equation}
\mathbb{P}(X \ge x)\ge\mathbb{P}(Z \ge x+\epsilon)-\mathbb{P}\bigl(|Y| \ge \epsilon\bigr).
\end{equation}
Moreover, from \eqref{E:nonuniform221},
%
%e3.24 #&#
\begin{equation}
\label{E:expectationesti} \mathbb{P}\bigl(|Y| \ge\epsilon\bigr)\le \frac{Cp_{\mathrm{2nd}}^m\sqrt{n}}{\epsilon\sqrt {kp_{\mathrm{max}}^m}},
\end{equation}
and from Lemma~\ref{lemma2},
\[
\mathbb{P}(Z \ge x\pm\epsilon)=\exp \bigl(-k\bigl(I_1(x\pm\epsilon
)+\mathrm{o}(1)\bigr) \bigr).
\]
Under the condition \eqref{E:intheorem},
%
%e3.25 #&#
\begin{equation}
\label{E:nonuniformratio} \frac{\mathbb{P}(|Y| \ge\epsilon)}{\mathbb{P}(Z \ge x\pm\epsilon)}\le \frac{Cp_{\mathrm{2nd}}^m\sqrt{n}}{\epsilon\sqrt{kp_{\mathrm{max}}^m}}\exp \bigl(k
\bigl(I_1(x\pm\epsilon)+\mathrm{o}(1)\bigr) \bigr) \to0, \qquad \mbox{as }
n \to\infty.
\end{equation}
Letting $\epsilon$ go to $0$, and repeating the arguments of the proof of
Theorem~\ref{T:diagramsmulti}, establishes \eqref{E:nonuniformresult1},
for any $x\ge2$,
under the conditions given in Theorem~\ref{T:nonuniform}.
For $\eqref{E:nonuniformresult0}$, for any $x < 2$ and $0 < \epsilon< 2-x$,
\[
\mathbb{P}(X\le x)\le\mathbb{P}(Z\le x+\epsilon)+\mathbb{P}\bigl(|Y|\ge \epsilon\bigr).
\]
From \eqref{E:nonuniformlemma1}, $\mathbb{P}(Z\le x+\epsilon)$ is exponentially
negligible with speed $k(m)$, and from arguments as in getting
\eqref{E:nonuniformratio}, $\mathbb{P}(|Y|\ge\epsilon)$ is also exponentially
negligible with speed $k(m)$, which proves \eqref
{E:nonuniformresult0}.
\end{pf*}
\begin{pf*}{Proof of Theorem~\ref{P:nonuniform} and Remark~\ref{remark1}}
Again, below, $m$ is short for $m(n)$ and $k$ is short for $k(m(n))$.
First, \eqref{E:nonuniformresult222} is a direct consequence of
\eqref{E:nonuniformresult1}, so let us prove \eqref{E:nonuniformresult2}.
When on the left of its simultaneous asymptotic mean,
$V_{1}^{\prime}(n,m)$ can be approximated by $F(n,k)$
(see \eqref{E:newbrownian}).
Hence, the rate function $K_{\eta}$ should be the corresponding ``left''
rate function of
the Brownian functional $F(1,k)$ (see \eqref{E:neweq}) with speed $k^2$.
From the right-hand side of \eqref{E:nonuniform22}, it is clear that
this new rate function
depends on $\eta= \lim_{n\to\infty}kp_{\mathrm{max}}^m$; let us denote it by
$K_{\eta}$.
Moreover, for $F(1,k)$, and from \cite{HX} or \cite{BGH},\vspace*{1pt}
\[
\sqrt{1-p_{\mathrm{max}}^m}F(1,k)\stackrel{\mathcal{L}}
{=}\tilde{\lambda}_1^0,
\]
where $\tilde{\lambda}_1^0$ is the largest eigenvalue of the diagonal block
corresponding to $p_{\mathrm{max}}^m$ in $\mathbf{X^0}$, and where $\mathbf
{X^0}$ is an
element of $\mathcal{G}^0(p_1^m,p_2^m,\ldots,p_m^m)$ (see \eqref{gene}).
So the rate function $K_{\eta}$ should also be the corresponding
``left'' rate
function for $\tilde{\lambda}_1^0$ with speed $k^2$.

Again, from \cite{HX},\vspace*{1pt}
%
%e3.26 #&#
\begin{equation}
\label{E:equivalence} \lambda_1^k\stackrel{\mathcal{L}} {=}
\tilde{\lambda}_1^0+\sqrt{p_{\mathrm{max}}^m}g,
\end{equation}
where $\lambda_1^k$ is the largest eigenvalue of an element of the
$k\times k$ GUE and
where $g$ is a standard normal random variable which is
independent of $\tilde{\lambda}_1^0$.

Let\vspace*{1pt}
%
%e3.27 #&#
%e3.28 #&#
\begin{eqnarray}
J(x)&=& %
\cases{ \displaystyle \inf_{\mu\in\mathcal{M}((-\infty,x])}I(\mu), &\quad \mbox{if $x\in
(-\infty ,2]$,}
\cr
0, &\quad \mbox{if $x\in[2,+\infty)$,} } %
\\
G_{\eta}(x)&=& %
\cases{ \displaystyle \frac{x^2}{2\eta}, &\quad \mbox{if $x\in
(-\infty,0]$,}
\cr
0, &\quad \mbox{if $x\in[0,+\infty)$,} } %
\end{eqnarray}
be the respective rate function for $\lambda_1^k$,
with speed $k^2$ and with $I(\mu)$ given in \eqref{traceless:ratefunction},
and for the Gaussian term.
Now, see \cite{DM}, when $x \le2$,\vspace*{1pt}
%
%e3.29 #&#
%e3.30 #&#
\begin{eqnarray}
J(x)&=&\tfrac{1}{216} \bigl(-x \bigl(-72x+x^3+30
\sqrt{12+x^2}+x^2\sqrt {12+x^2} \bigr)\nonumber
\\[-8pt]\\[-8pt]
&&\hphantom{\tfrac{1}{216} \bigl(}{}-216\log \bigl(\tfrac{1}{6} \bigl(x+\sqrt{12+x^2} \bigr)
\bigr) \bigr).\nonumber
\end{eqnarray}
Hence,
%
%e3.31 #&#
\begin{eqnarray}\label{ratefunction_second}
J^\prime(x)&=&\tfrac{1}{54} \bigl(-x^3+36x-
\bigl(12+x^2\bigr)^{3/2} \bigr),
\\
 J^{\prime\prime}(x)&=&\tfrac{1}{18} \bigl(12-x^2-x
\sqrt{12+x^2} \bigr).
\end{eqnarray}
Clearly, for $x \in(-\infty,2)$, $0<J^{\prime\prime}(x)<1$ and by a
Taylor expansions
for $J$ and $J^\prime$, and for $x<-5$,
%
%e3.33 #&#
\begin{eqnarray}
\label{expan-rate} J(x)&=&\frac{x^2}{2}+\log(-x)+\frac{3}{4}+e_1(x),
\\
\label{expan-derivative} J^\prime(x)&=&x+\frac{1}{x}+e_2(x),
\end{eqnarray}
with $|e_1(x)|\le2/x^2$ and $|e_2(x)|\le4/|x|^3$.

From \eqref{E:equivalence}, it is well known (see \cite{DemboZ,Ro}) that,
%
%e3.35 #&#
\begin{equation}
\label{inf-con} J(x)=K_{\eta} \,\Box\, G_{\eta}(x):=\inf
_{y \in\mathbb{R}}\bigl(K_{\eta
}(y)+G_{\eta}(x-y)\bigr),
\end{equation}
and taking Legendre transforms:
\[
K_{\eta}(x)= \bigl(J^*(y)-G_{\eta}^*(y) \bigr)^*(x),
\]
where
\[
G^*(y)= %
\cases{ \displaystyle \frac{\eta y^2}{2}, &\quad \mbox{if $y \le0$,}
\cr
+
\infty, &\quad \mbox{if $y>0$,} } %
\]
so that
%
%e3.36 #&#
\begin{equation}
\label{ratewithbeta} K_{\eta}(x)=\sup_{y\le0}
\biggl(xy-J^*(y)+\frac{\eta y^2}{2} \biggr).
\end{equation}
Therefore, for $0< \eta<1$, $K_{\eta}$ interpolates between $K_0=J$
and $K_1=K$.
From the very definition of the Legendre transform,
\[
J^*(y)=\sup_{x\in\mathbb{R}} \bigl(xy-J(x) \bigr),
\]
there exists, for each $y\le0$, a unique solution, denoted by $S(y)$,
to $J^\prime(x)=y$ for $x\in(-\infty,2]$. Clearly, the function
$S$ is increasing on $(-\infty,0]$, with $S(0)=2$, $\lim_{y \to-\infty
}S(y)=-\infty$,
and with
\[
S^\prime(y)=\frac{1}{J^{\prime\prime}(S(y))},
\]
for $y<0$. Thus, for $y\le2$,
\[
J^*(y)=yS(y)-J\bigl(S(y)\bigr),
\]
and,
\[
K_{\eta}(x)=\sup_{y\le0} \biggl(xy-yS(y)+J\bigl(S(y)
\bigr)+\frac{\eta y^2}{2} \biggr).
\]
For $y\le0$, let
\[
H_{x,\eta}(y):=xy-yS(y)+J\bigl(S(y)\bigr)+\frac{\eta y^2}{2},
\]
then
\[
H_{x,\eta}^\prime(y)=x-S(y)+\eta y,\qquad  H_{x,\eta}^{\prime\prime}(y)
=-\frac{1}{J^{\prime\prime}(S(y))}+\eta,
\]
so $H_{x,\eta}^{\prime\prime}(y)<0$ for $y\in(-\infty,0)$, $x\in
\mathbb{R}$
and $0\le\eta\le1$. When $x\ge2$, for any $0\le\eta\le1$,
$H_{x,\eta}^\prime(y)>0$ for
$y<0$ with $H_{x,\eta}^\prime(0)\ge0$, thus $K_{\eta}(x)=\sup_{y\le
0}H_{x,\eta}(y)=H_{x,\eta}(0)=0$.
Let us now deal with $x<2$. First, from \eqref{expan-derivative}, it
can be
shown that for $y<-6$,
\[
y< S(y)< y+1,
\]
and thus since $x-J^\prime(x)$ is increasing on $(-\infty,2]$,
\[
S(y)-y=S(y)-J^\prime\bigl(S(y)\bigr)<y+1-J^\prime(y+1)<-
\frac{2}{y+1},
\]
which further yields
\[
y<S(y)<y-\frac{2}{y+1}.
\]
Moreover, when $y<-6$,
%
%e3.37 #&#
\begin{eqnarray}
\label{rate-estimate} \biggl\llvert H_{x,\eta}(y)- \biggl(xy-y^2+J(y)+
\frac{\eta y^2}{2} \biggr) \biggr\rrvert &\le&\llvert y \rrvert \bigl\llvert S(y)-y
\bigr\rrvert + \bigl\llvert J\bigl(S(y)\bigr)-J(y)\bigr\rrvert
\nonumber
\\
&\le&2\biggl\llvert \frac{y}{y+1}\biggr\rrvert +\bigl\llvert
J^\prime(y+1)\bigr\rrvert \bigl\llvert S(y)-y\bigr\rrvert
\\
&\le&3+3=6.\nonumber
\end{eqnarray}
Combining \eqref{rate-estimate} with \eqref{expan-rate}, it follows
that for $y<-6$,
%
%e3.38 #&#
\begin{equation}
\label{rate-inequality} \biggl\llvert H_{x,\eta}(y)- \biggl(xy+\log(-y)-
\frac{1-\eta}{2}y^2 \biggr)\biggr\rrvert \le7.
\end{equation}
When $\eta=1$, for any $x\le0$, $H_{x,1}^\prime(y)< 0$ for $y\le0$, thus
\[
K_{1}(x)=\lim_{y\to-\infty}H_{x,1}(y)=+\infty.
\]
For $0<x<2$, since $S(y)-y$ is increasing on $(-\infty,0]$ with a
range of $(0,2]$, there exists a unique solution, denoted by $T_1(x)$,
to $H_{x,1}^\prime(y)=x-S(y)+y=0$.
Note that\vadjust{\goodbreak} $y=T_1(x)$ is the maximizer of $H_{x,1}(y)$ and as
$x \to0$, $T_1(x) \to-\infty$, thus there exists $\delta>0$, such
that for $x<\delta$,
\[
K_{1}(x)=\sup_{y\le-6}H_{x,1}(y).
\]
Since for $x<1/6$,
\[
\sup_{y\le-6} \bigl(xy+\log(-y) \bigr)=-1-\log x,
\]
when combined with \eqref{rate-inequality}, it follows that
for $x$ close enough to 0,
\[
\bigl\llvert K_1(x)-(-\log x)\bigr\rrvert \le8.
\]

When $0<\eta<1$, for any $x<2$, there exists a unique
solution, denoted by $T_{\eta}(x)$, to $H_{x,\eta}^\prime
(y)=x-S(y)+\eta y=0$.
Again, $y=T_{\eta}(x)$ is the maximizer of $H_{x,\eta}(y)$ and as
$x \to-\infty$, $T_{\eta}(x) \to-\infty$. Repeating the arguments of
the case $\eta=1$, gives as $x \to-\infty$,
\[
K_{\eta}(x)\sim\frac{x^2}{2 (1-\eta)}+\log \biggl(-\frac{x}{1-\eta}
\biggr).
\]
This last statement (clearly consistent with the case $\eta=0$),
finishes to prove the last assertions of Remark~\ref{remark1}.
The rest of the proof of Theorem~\ref{P:nonuniform} follows along the
lines of
the proofs of Lemma~\ref{lemma2} and of Theorem~\ref{T:nonuniform}, and
is therefore left to the interested reader.
\end{pf*}

%s4 #&#
\section{Proof of Theorem \texorpdfstring{\protect\ref{T:concentration}}{1.5} 
and Theorem \texorpdfstring{\protect\ref
{T:nonuniformconcentration}}{1.6}}
Left and right concentration inequalities for the largest eigenvalue
$\lambda_1^m$ of an element of the $m\times m$ GUE are respectively given
in \cite{Au} and \cite{LR}. More precisely:
%
%pr4.1 #&#
\begin{proposition}\label{P:concentration}
Let $m\ge1$ and let $\epsilon>0$, then for some absolute constant $C_0>0$,
%
%e4.1 #&#
\begin{equation}
\label{E:concentrationpropo1} \mathbb{P}\bigl(\lambda_1^m\ge2
\sqrt{m}(1+\epsilon)\bigr)\le C_0 \mathrm{e}^{-m\epsilon
^{3/2}/C_0}.
\end{equation}
Likewise, for some absolute constant $\bar{C}_0>0$, and all $m\ge1$
and $0<\epsilon\le1$,
%
%e4.2 #&#
\begin{equation}
\label{E:concentrationpropo2} \mathbb{P}\bigl(\lambda_1^m\le2
\sqrt{m}(1-\epsilon)\bigr)\le\bar{C}_0 \mathrm{e}^{-m^2\epsilon^{3}/\bar{C}_0}.
\end{equation}
\end{proposition}
Next, to prove \eqref{E:concentrationthm1}, assume first
that $m\epsilon^{3/2}\ge1$. Then for any $0<\epsilon<1$,
%
%e4.3 #&#
\begin{eqnarray}\label{E:concentration1}
&&\mathbb{P} \biggl(\frac{V_1(n,m)-n/m}{\sqrt{n/m}} \ge2\sqrt {m}(1+\epsilon) \biggr)
\nonumber
\\
&&\quad  \le\mathbb{P} \biggl(\sqrt{\frac{m-1}{m}}\frac{\tilde{L}_1(n,m)}{2\sqrt{mn}} \ge1+
\frac{\epsilon}{2} \biggr)
\\
&&\qquad {} +\mathbb{P} \biggl(\sqrt{\frac{m-1}{m}}\frac{|\tilde{V}_{1}(n,m)-\tilde
{L}_1(n,m)|}{2\sqrt{mn}}\ge
\frac{\epsilon}{2} \biggr).
\nonumber
\end{eqnarray}
As before,
\[
\sqrt{\frac{m-1}{m}}\frac{\tilde{L}_1(n,m)}{\sqrt{n}}\stackrel{\mathcal {L}} {=}
\lambda_1^{m,0}
\]
and
\[
\lambda_{1}^{m}\stackrel{\mathcal{L}} {=}
\lambda_{1}^{m,0}+g_{m},
\]
where $g_{m}$ is a centered Gaussian random variable with variance $1/m$,
independent of $\lambda_{1}^{m,0}$. So,
\begin{eqnarray*}
\mathbb{P} \biggl( \sqrt{\frac{m-1}{m}}\frac{\tilde{L}_1(n,m)}{2\sqrt {mn}}\ge1+
\frac{\epsilon}{2} \biggr)&\le&\mathbb{P} \biggl(\lambda _1^{m}
\ge2\sqrt{m} \biggl(1+\frac{\epsilon}{4} \biggr) \biggr)+\mathbb{P} \biggl(
g_m\ge\frac{\sqrt{m}\epsilon}{2} \biggr)
\\
&\le& C_1 \mathrm{e}^{-m\epsilon^{3/2}/C_1}+C_1
\mathrm{e}^{-m^2\epsilon^{2}/C_1},
\end{eqnarray*}
for some positive constant $C_1$. Now from \eqref{E:concentrationback},
\eqref{E:diagramsuniform1} and
\eqref{E:diagramsuniform2}, the second term on the right-hand side of
\eqref{E:concentration1} is upper-bounded by:
\[
\mathbb{P} \biggl(\frac{|\tilde{V}_{1}(n,m)-\tilde{L}_1(n,m)|}{2\sqrt {mn}}\ge\frac{\epsilon}{2} \biggr) \le
C_2\sqrt{mn}\mathrm{e}^{-\sqrt{n}\epsilon/C_2m}+ C_2mn\mathrm{e}^{-n\epsilon^2/C_2m}.
\]
In order to reach \eqref{E:concentrationthm1}, we need to show that there
exists a positive constant $C(A,\alpha)$, depending only on
$A$ and $\alpha$, such that
%
%e4.4 #&#
%e4.5 #&#
%e4.6 #&#
\begin{eqnarray}
\label{E:concentration11}C(A,\alpha) \mathrm{e}^{-m\epsilon^{3/2}/C(A,\alpha)}&\ge& C_1 \mathrm{e}^{-m^2\epsilon
^{2}/C_1},
\\
\label{E:concentration12}C(A,\alpha) \mathrm{e}^{-m\epsilon^{3/2}/C(A,\alpha)}&\ge& C_2\sqrt{mn}
\mathrm{e}^{-\sqrt
{n}\epsilon/C_2m},
\\
\label{E:concentration13}C(A,\alpha) \mathrm{e}^{-m\epsilon^{3/2}/C(A,\alpha)}&\ge& C_2mn\mathrm{e}^{-n\epsilon
^2/C_2m}.
\end{eqnarray}
First, since $m\epsilon^{3/2}\ge1$, \eqref{E:concentration11} is
satisfied by
choosing $C(A,\alpha) \ge C_1$.
Now taking logarithms in \eqref{E:concentration12},
$C(A,\alpha)$ has to be such that:
%
%e4.7 #&#
\begin{equation}
\log\frac{C_2}{C(A,\alpha)}+\frac{1}{2}\log( mn)\le m \epsilon
^{3/2} \biggl(-\frac{1}{C(A,\alpha)}+\frac{\sqrt{n}}{C_2 m^2 \epsilon
^{1/2}} \biggr).
\end{equation}
Moreover, under the condition $m\le An^{\alpha}$, we have:
\[
\frac{\sqrt{n}}{C_2 m^2 \epsilon^{1/2}}\ge \frac{\sqrt{n}}{C_2 m^2 }\ge\frac{n^{\sklfrac{1}{2}-2\alpha}}{A^2C_2}.
\]
Therefore, if $\alpha<1/4$, it is enough to choose $C(A,\alpha)$ satisfying
\[
\log\frac{\sqrt{A}C_2}{C(A,\alpha)}+\frac{1}{C(A,\alpha)} \le\frac{n^{\sklfrac{1}{2}-2\alpha}}{A^2C_2}-
\frac{1+\alpha}{2}\log n.
\]
Since for all integers $n\ge1$,
\[
\frac{n^{\sfrac{1}{2}-2\alpha}}{A^2C_2}-\frac{1+\alpha}{2}\log n \ge\frac{1+\alpha}{1-4\alpha} \biggl(1-
\log\frac{A^2C_2 (1+\alpha)}{1-4
\alpha} \biggr),
\]
we just need to guarantee that
%
%e4.8 #&#
\begin{equation}
\label{con:1} \log\frac{\sqrt{A}C_2}{C(A,\alpha)}+\frac{1}{C(A,\alpha)} \le\frac{1+\alpha}{1-4\alpha}
\biggl(1-\log\frac{A^2C_2 (1+\alpha)}{1-4
\alpha} \biggr).
\end{equation}
But, from our choice of $\alpha$, $(1+\alpha)/(1-4\alpha)>1$, so by choosing
%
%e4.9 #&#
\begin{equation}
\label{con:4} C(A,\alpha)\ge C\max\bigl(A^{5/2},1\bigr)
\frac{1+\alpha}{1-4\alpha} \exp \biggl(\frac{1+\alpha}{1-4\alpha} \biggr),
\end{equation}
for some large enough absolute constant $C>0$,
\eqref{con:1} and \eqref{E:concentration12} are satisfied.
Finally, by taking logarithms, \eqref{E:concentration13} becomes,
%
%e4.10 #&#
\begin{equation}
\label{E:conesti1} \log\frac{C_2}{C(A,\alpha)}+\log(mn)\le m \epsilon^{3/2}
\biggl(-\frac
{1}{C(A,\alpha)}+\frac{n \epsilon^{1/2}}{C_2 m^2} \biggr).
\end{equation}
From the condition $m\le An^{\alpha}$, we just need,
%
%e4.11 #&#
\begin{equation}
\label{con:2} \log\frac{AC_2}{C(A,\alpha)}+\frac{1}{C(A,\alpha)}\le\frac
{1}{A^{7/3}C_2}
n^{1-\sfrac{7\alpha}{3}}-(1+\alpha)\log n.
\end{equation}
Now repeating the previous arguments, taking the minimum on the right-hand
side of \eqref{con:2}, it follows that
%
%e4.12 #&#
\begin{equation}
\label{con:3} \log\frac{AC_2}{C(A,\alpha)}+\frac{1}{C(A,\alpha)} \le
\frac{1+\alpha}{1-7\alpha/3} \biggl(1-\log\frac{A^{7/3}C_2(1+\alpha
)}{1-7\alpha/3} \biggr).
\end{equation}
Again, for $0< \alpha< 1/4$, $1 < (1+\alpha)/(1-7\alpha/3) < 3$, so as
long as
%
%e4.13 #&#
\begin{equation}
\label{con:6} C(A,\alpha)\ge C\max\bigl(A^{10/3},1\bigr)
\frac{1+\alpha}{1-7\alpha/3} \exp \biggl(\frac{1+\alpha}{1-7\alpha/3} \biggr),
\end{equation}
for some large enough absolute constant $C$, then $C(A,\alpha)$ will
also satisfy \eqref{con:3} and therefore also \eqref{E:concentration13}.

Combining \eqref{con:4} and \eqref{con:6},
if $m\epsilon^{3/2}\ge1$, and $m\le An^{\alpha}$,
with $\alpha<1/4$, there exist a positive constant
%
%e4.14 #&#
\begin{equation}
\label{con:7} C(A,\alpha)= C\max\bigl(A^{10/3},1\bigr)
\frac{1+\alpha}{1-4\alpha} \exp \biggl(\frac{1+\alpha}{1-4\alpha} \biggr),
\end{equation}
so that \eqref{E:concentrationthm1} holds true
for all $0 < \epsilon< 1$. When $m\epsilon^{3/2}< 1$,
\[
C(A,\alpha) \mathrm{e}^{-m\epsilon^{3/2}/C(A,\alpha)}\ge C \mathrm{e}^{-1/C}\ge1,
\]
as $C$ is large enough, and \eqref{E:concentrationthm1} follows.
So combining these two cases, there exists $C(A,\alpha)$ as in \eqref
{con:7}, with
$C$ large enough, such that \eqref{E:concentrationthm1} is satisfied.

Likewise, for the proof of \eqref{E:concentrationthm2},
first assume that $m^2\epsilon^3\ge1$, and\vspace*{-1.5pt}
%
%e4.15 #&#
\begin{eqnarray}
\label{E:concentration2} &&\mathbb{P} \biggl(\frac{V_1(n,m)-n/m}{\sqrt{n/m}} \ge2\sqrt {m}(1-\epsilon)
\biggr) \nonumber
\\
&&\quad \le\mathbb{P} \biggl(\sqrt{\frac{m-1}{m}}\frac{\tilde{L}_1(n,m)}{2\sqrt {mn}}\le1-
\frac{\epsilon}{2} \biggr)
\nonumber
\\[-8pt]\\[-8pt]
&&\qquad {} +\mathbb{P} \biggl(\sqrt{\frac{m-1}{m}}\frac{|\tilde{V}_{1}(n,m)-\tilde
{L}_1(n,m)|} {
2\sqrt{mn}}\ge
\frac{\epsilon}{2} \biggr)
\nonumber
\\
&&\quad \le C_1 \mathrm{e}^{-m^2\epsilon^{3}/C_1}+C_1
\mathrm{e}^{-m^2\epsilon^{2}/C_1}+C_2\sqrt {mn}\mathrm{e}^{-\sqrt{n}\epsilon/C_2m}+
C_2mn\mathrm{e}^{-n\epsilon^2/C_2m}.
\nonumber
\end{eqnarray}
Repeating previous arguments, and as long as $m\le An^{\alpha}$, with
$\alpha<1/6$,
there exists a positive constant\vspace*{-1.5pt}
\[
\bar{C}(A,\alpha)=\bar{C}\max\bigl(A^4,1\bigr) \frac{1+\alpha}{1-6\alpha}
\exp \biggl(\frac{1+\alpha}{1-6\alpha} \biggr),
\]
so that \eqref{E:concentrationthm2} is satisfied.
Once more, taking $\bar{C}$ large enough, the case $m^2\epsilon^3< 1$ follows,
and \eqref{E:concentrationthm2} is proved.

The proof for the non-uniform case is similar to the uniform one.
For \eqref{E:nonuniformconcentrationthm1}, assume at first that
$k\epsilon^{3/2}\ge1$, then\vspace*{-1.5pt}
\begin{eqnarray*}
&&\mathbb{P} \biggl(\frac{V_1(n,m)-np_{\mathrm{max}}^m}{\sqrt{nkp_{\mathrm{max}}^m}}\ge 2(1+\epsilon) \biggr)
\\
&&\quad  \le\mathbb{P} \biggl( \frac{V_1(n,m)-V_1^{\prime}(n,m)}{2\sqrt {nkp_{\mathrm{max}}^m}}\ge\frac{\epsilon}{3} \biggr)+
\mathbb{P} \biggl( \sqrt {1-p_{\mathrm{max}}^m}
\frac{\tilde{V}_{1}^{\prime}(n,m)-F(n,k)}{2\sqrt{nk}}\ge \frac{\epsilon}{3} \biggr)
\nonumber
\\
&&\qquad {} + \mathbb{P} \biggl(\sqrt{1-p_{\mathrm{max}}^m}
\frac{F(n,k)}{2\sqrt{nk}}\ge1+\frac
{\epsilon}{3} \biggr)
\\
&&\quad  =A_1+A_2+A_3.
\end{eqnarray*}
From \eqref{E:expectationesti}, \eqref{E:nonuniform-concentrationback}
and \eqref{E:nonuniform22},\vspace*{-1.5pt}
\begin{eqnarray*}
A_1 &\le&\frac{C_1p_{\mathrm{2nd}}^m\sqrt{n}}{\epsilon\sqrt{kp_{\mathrm{max}}^m}},
\nonumber
\\
A_2 &\le& C_2 nk \exp \biggl(-\frac{n\epsilon^2}{C_2 k}
\biggr)+C_2 k\sqrt {np_{\mathrm{max}}^m}\exp
\biggl(-\frac{\sqrt{np_{\mathrm{max}}^m}\epsilon}{C_2\sqrt {k}} \biggr),
\nonumber
\\
A_3 &\le&\mathbb{P} \biggl( Z_{k}\ge\frac{\epsilon}{3}
\biggr)+\mathbb {P} \biggl( \lambda_1^{k}\ge2 \biggl(1+
\frac{\epsilon}{6} \biggr) \biggr)
\nonumber
\\
&\le& C_3 \exp \biggl(-\frac{k^2\epsilon^2}{C_3} \biggr) +C_3
\exp \biggl(-\frac{k\epsilon^{3/2}}{C_3} \biggr).
\nonumber
\end{eqnarray*}
In order to reach \eqref{E:nonuniformconcentrationthm1}, we need to
show that there exists
a positive constant $C(A,B,\alpha)$, depending only on $A$, $B$ and
$\alpha$, such that
%
%e4.16 #&#
%e4.17 #&#
%e4.18 #&#
%e4.19 #&#
\begin{eqnarray}
\label{E:nonuniformconcentration1}C(A,B,\alpha) \exp \biggl(-\frac{k\epsilon^{3/2}}{C(A,B,\alpha)} \biggr)&\ge&\frac{C_1p_{\mathrm{2nd}}^m\sqrt{n}}{\epsilon\sqrt {kp_{\mathrm{max}}^m}},
\\
\label
{E:nonuniformconcentration2}C(A,B,\alpha) \exp \biggl(-\frac{k\epsilon^{3/2}}{C(A,B,\alpha)} \biggr)&\ge& C_2 nk
\exp \biggl(-\frac{n\epsilon^2}{C_2 k} \biggr),
\\
\label{E:nonuniformconcentration3}C(A,B,\alpha) \exp \biggl(-\frac{k\epsilon^{3/2}}{C(A,B,\alpha)} \biggr)&\ge& C_2 k
\sqrt{np_{\mathrm{max}}^m}\exp \biggl(-\frac{\sqrt{np_{\mathrm{max}}^m}\epsilon
}{C_2\sqrt{k}}
\biggr),
\\
\label
{E:nonuniformconcentration4}C(A,B,\alpha) \exp \biggl(-\frac{k\epsilon^{3/2}}{C(A,B,\alpha)} \biggr)&\ge &C_3 \exp
\biggl(-\frac{k^2\epsilon^2}{C_3} \biggr).
\end{eqnarray}

First, taking logarithms in \eqref{E:nonuniformconcentration3}, gives:
\[
\log\frac{C_2}{C(A,B,\alpha)}+\log k +\frac{1}{2}\log\bigl(np_{\mathrm{max}}^m
\bigr)\le k \epsilon^{3/2} \biggl(-\frac{1}{C(A,B,\alpha)}+
\frac{\sqrt {np_{\mathrm{max}}^m}}{C_2\sqrt{\epsilon k^3}} \biggr).
\]
Next,
\[
\frac{\sqrt{np_{\mathrm{max}}^m}}{C_2\sqrt{\epsilon k^3}}\ge \frac{\sqrt{(np_{\mathrm{max}}^m)^{1-3/\alpha}}}{A^{3/2\alpha}C_2},
\]
so if $\alpha>3$, then there exists a constant $C(A,B,\alpha)>0$,
satisfying \eqref{E:nonuniformconcentration3}.
In fact, here $C(A,B,\alpha)$ just needs to be such that
\[
\log\frac{A^{1/\alpha}C_2}{C(A,B,\alpha)}+\frac{1} {
C(A,B,\alpha)}\le\frac{\alpha+2}{\alpha-3} \biggl(1-\log
\frac{A^{3/2\alpha}C_2(\alpha+2)}{\alpha-3} \biggr),
\]
which forces
%
%e4.20 #&#
\begin{equation}
\label{con:5} C(A,B,\alpha)\ge C \max\bigl(A^{2/\alpha},1\bigr)
\frac{\alpha+2}{\alpha-3} \exp \biggl(\frac{\alpha+2}{\alpha-3} \biggr),
\end{equation}
for a large enough absolute constant $C>0$.

Second, taking logarithms in \eqref{E:nonuniformconcentration1}, gives:
\[
\log\frac{C_1}{C(A,B,\alpha)} +\log \biggl(\frac{p_{\mathrm{2nd}}^m\sqrt{n}}{\sqrt{kp_{\mathrm{max}}^m}} \biggr)\le -
\frac{k \epsilon^{3/2}}{C(A,B,\alpha)}+\log\epsilon.
\]
From \eqref{E:nonconexp1} and the assumption $k\epsilon^{3/2}\ge1$, in
order for \eqref{E:nonuniformconcentration1} to hold true, $C(A,B,\alpha
)$ needs to satisfy
\[
\log\frac{C_1\sqrt{B}}{C(A,B,\alpha)}-\frac{k}{2} \le-\frac
{k}{C(A,B,\alpha)}-
\frac{2}{3}\log k,
\]
which further forces
%
%e4.21 #&#
\begin{equation}
\label{con:9} C(A,B,\alpha)\ge C \max(\sqrt{B},1),
\end{equation}
with the absolute constant $C$ large enough.

For \eqref{E:nonuniformconcentration2}, as we did with \eqref
{E:concentration13}, and under the condition $k^{\alpha}/p_{\mathrm{max}}^m\le
An$ with $\alpha>3$, we need:
%
%e4.22 #&#
\begin{equation}
\label{con:8} C(A,B,\alpha)\ge C\max\bigl(A^{10/3\alpha},1\bigr)
\frac{3\alpha+3}{3\alpha-7} \exp \biggl(\frac{3\alpha+3}{3\alpha-7} \biggr),
\end{equation}
with the absolute constant $C$ large enough.
Finally, \eqref{E:nonuniformconcentration4} is easily satisfied
since $k \epsilon^{3/2}\ge1$.
Moreover, when $k \epsilon^{3/2}< 1$,
then \eqref{E:nonuniformconcentrationthm1} holds, given $C>0$ large enough.
Combining \eqref{con:5}, \eqref{con:9} and \eqref{con:8}, choosing
\[
C(A,B,\alpha)= C \max\bigl(A^{10/3\alpha},1\bigr) \max(\sqrt{B},1)
\frac{\alpha+2}{\alpha-3} \exp \biggl(\frac{\alpha
+2}{\alpha-3} \biggr),
\]
with $C>0$, some large enough absolute constant, then
\eqref{E:nonuniformconcentrationthm1} holds under the given conditions.
Likewise, we can prove \eqref{E:nonuniformconcentrationthm2}.

\begin{appendix}
%s5 #&#
\section*{Appendix: Large deviations for the spectrum of the traceless GUE}\label{app}
\setcounter{equation}{0}
For any integer $m\ge2$, let the random matrix $\mathbf{X}$ be an
element of the $m \times m$ GUE. Let $(\lambda_{1},\lambda_{2},\ldots,
\lambda_{m})$ be the spectrum of $\mathbf{X}$, and let
\[
(\xi_{1},\xi_{2},\ldots, \xi_{m})=
\frac{1}{\sqrt{m}}(\lambda_{1},\lambda _{2},\ldots,
\lambda_{m}).
\]
The joint probability density of $(\xi_{1},\xi_{2},\ldots, \xi_{m})$ is
given by
%
%e5.1 #&#
\begin{equation}
\label{density} \phi_{m}(\xi_{1},\xi_{2},\ldots,
\xi_{m})=\frac{1}{Z_m}\exp \Biggl(-\frac
{m}{2}\sum
_{i=1}^m \xi_i^2
\Biggr)\prod_{1\le i<j \le m}(\xi_i-
\xi_j)^2,
\end{equation}
where
%
%e5.2 #&#
\begin{equation}
Z_m=(2\uppi)^{\sfrac{m}{2}}m^{-\sfrac{m^2}{2}}\prod
_{j=1}^m j!,
\end{equation}
see Theorem~2.5.2 in \cite{AGZ} and also Theorem~3.3.1 in \cite{Me}.

Let $(\lambda_{1}^{0},\lambda_{2}^{0},\ldots, \lambda_{m}^{0})$ be
the spectrum of $\mathbf{X}-\tr(\mathbf{X})/m$, an element of the $m
\times m$
traceless GUE, and again, let
\[
\bigl(\xi_{1}^{0},\xi_{2}^{0},
\ldots, \xi_{m}^{0}\bigr)=\frac{1}{\sqrt{m}}\bigl(\lambda
_{1}^{0},\lambda_{2}^{0},\ldots,
\lambda_{m}^{0}\bigr).
\]
The joint distribution function of $(\xi_{1}^{0},\xi_{2}^{0},\ldots, \xi
_{m}^{0})$ is
given by
%
%e5.3 #&#
\begin{eqnarray}
&&\mathbb{P} \bigl(\xi_1^0\le s_1,
\xi_2^0\le s_2,\ldots,\xi_m^0
\le s_m \bigr)\nonumber
\\[-8pt]\\[-8pt]
&&\quad  =\sqrt{2\uppi}\int_{\mathcal{L}(s_1,\ldots,s_m)} \phi_{m}(x_1,x_2,
\ldots,x_m)\,\mathrm{d}x_1\cdots \,\mathrm{d}x_{m-1},
\nonumber
\end{eqnarray}
where
%
%e5.4 #&#
\begin{eqnarray*}
\mathcal{L}(s_1,\ldots,s_m)&:=& \Biggl
\{x=(x_1,\ldots,x_m)\in\mathbb {R}^m\dvtx \sum
_{i=1}^m x_i=0, \mbox{and }
x_i \le s_i,
\\
&&\hphantom{\Biggl
\{}\mbox{for each }i=1,\ldots,m \Biggr\}.
\nonumber
\end{eqnarray*}

Let $(\xi_{1}^{m},\xi_{2}^{m},\ldots, \xi_{m}^{m})$ be
the non-increasing rearrangement of $(\xi_{1},\xi_{2},\ldots, \xi_{m})$,
and let $(\xi_{1}^{m,0},\allowbreak \xi_{2}^{m,0},\ldots, \xi_{m}^{m,0})$ be
the non-increasing rearrangement of
$(\xi_{1}^{0},\xi_{2}^{0},\ldots, \xi_{m}^{0})$, then, for
example, see~\cite{HX},
%
%e5.5 #&#
\begin{equation}
\label{E:equaleigen} \bigl(\xi_{1}^{m},\xi_{2}^{m},
\ldots, \xi_{m}^{m}\bigr)\stackrel{\mathcal {L}} {=}\bigl(
\xi_{1}^{m,0},\xi_{2}^{m,0},\ldots,
\xi_{m}^{m,0}\bigr)+g_{m}\mathbf{e}_{m},
\end{equation}
where $g_m$ is a centered Gaussian random variable
with variance $1/m^2$, independent of the vector $(\xi_{1}^{m,0},\xi
_{2}^{m,0},\ldots, \xi_{m}^{m,0})$,
and where $\mathbf{e}_{m}=(1,1,\ldots,1)$.

As shown in \cite{BG}, the law of the spectral
measure $\hat{\mu}^m=\frac{1}{m}\sum_{i=1}^m\delta_{\xi_i}$ satisfies a large
deviation principle on the set
$\mathcal{P}(\mathbb{R})$ of probability measures on $\mathbb{R}$, and
with good rate function $I$, in the scale $m^2$. Moreover, $I$ is given by
%
%e5.6 #&#
\begin{equation}
\label{traceless:ratefunction} I(\mu)=\frac{1}{2}\int x^2 \mu(\mathrm{d}x)-
\iint\log|x-y|\mu(\mathrm{d}x)\mu(\mathrm{d}y)-\frac
{3}{4} ,
\end{equation}
and its unique minimizer is the semicircular probability measure
\[
\sigma(\mathrm{d}x)=\frac{1}{2\uppi}\mathbf{1}_{|x|\leq2}
\sqrt{4-x^2}\,\mathrm{d}x.
\]

Based on this LDP for $\hat{\mu}^m$, the LDP for the largest (or $r$th
largest) eigenvalue of
the GOE with an explicit rate function is obtained in \cite{BDG} and
\cite{ABC}
(see also \cite{Jo1} for generalizations).
Following the approach and the techniques developed there, and taking into
account \eqref{E:equaleigen}, we get a multidimensional LDP for the first
$r$ eigenvalues of the traceless GUE:
%
%th5.1 #&#
\begin{theorem}\label{T:multi1}
Let $r\in\mathbb{N}$, on $\mathcal{L}^r:=\{(x_1,x_2,\ldots,x_r)\in
\mathbb{R}^r\dvtx x_1\ge x_2 \ge\cdots\ge x_r\}$,
$(\xi_{1}^{m,0},\break \xi_{2}^{m,0},\ldots,  \xi_{r}^{m,0})$ satisfies a
LDP with speed $m$ and a good rate function
\[
I_r(x_1,x_2,\ldots,x_r)=
\cases{ 2\displaystyle \sum_{i=1}^r\displaystyle \int
_2^{x_i}\sqrt{(z/2)^2-1}\,\mathrm{d}z,
&\quad \mbox{if $x_1\ge x_2\ge \cdots\ge x_r \ge2$,}
\cr
+\infty, &\quad \mbox{otherwise.} } %
\]
\end{theorem}
\begin{pf} Let
\[
Q_{m}(\mathrm{d}\xi_{1},\mathrm{d}\xi_{2},
\ldots, \mathrm{d}\xi_{m})=\frac{1}{Z_m}\exp \Biggl(-
\frac
{m}{2}\sum_{i=1}^m
\xi_i^2 \Biggr)\prod_{1\le i<j \le m}(
\xi_i-\xi _j)^2\prod
_{i=1}^m \,\mathrm{d}\xi_i.
\]
From \cite{BDG}, $(\xi_{1}^{m},\xi_{2}^{m},\ldots, \xi_{r}^{m})$
satisfies a LDP with speed
$m$ and rate function $I_r$ on $\mathcal{L}^r$.
To prove the validity of the same results for
$(\xi_{1}^{m,0},\xi_{2}^{m,0},\ldots, \xi_{r}^{m,0})$, it is enough
to show that\vspace*{1pt}
%
%e5.7 #&#
\begin{equation}\label{multi2:infinity}
\limsup_{m \to\infty}\frac{1}{m}\log Q_m\bigl(
\xi_r^{m,0}\le x\bigr)=-\infty ,
\end{equation}
for any $x<2$, and since $I_r(x_1,x_2,\ldots,x_r)$ is continuous,
increasing in each individual
variable, on $\mathcal{L}^r\cap[2,\infty)^r$,
%
%e5.8 #&#
\begin{eqnarray}\label{multi2:finite}
\lim_{m \to\infty} \frac{1}{m} \log Q_m \bigl(
\xi_{1}^{m,0} \ge x_1,\ldots,
\xi_{r}^{m,0} \ge x_r \bigr)=-2\sum
_{i=1}^r \int_2^{x_i}
\sqrt{(z/2)^2-1}\,\mathrm{d}z,
\end{eqnarray}
for all $x_1\ge x_2\ge \cdots\ge x_r \ge2$.

First, for $x<2$, let $\delta=2-x$, so
%
%e5.9 #&#
\begin{eqnarray}
Q_m\bigl(\xi_r^{m,0}\le x\bigr)&\le&
Q_m\bigl(\xi_r^{m,0}+g_m\le x+
\delta/2\bigr)+\mathbb {P}(g_m\ge\delta/2)\nonumber
\\[-8pt]\\[-8pt]
&=&Q_m\bigl(\xi_r^{m}\le x+\delta/2\bigr)+
\mathbb{P}(g_m\ge\delta/2).
\nonumber
\end{eqnarray}
Since,
%
%e5.10 #&#
\begin{equation}\label{multi:traceless1}
\mathbb{P} (g_m\ge\delta )\sim\frac{1}{\sqrt{2 \uppi}m \delta
}
\mathrm{e}^{-m^2 \delta^2/2}, \qquad \mbox{as }m \to\infty,
\end{equation}
\eqref{multi2:infinity} follows.
For \eqref{multi2:finite},
fix $x_1\ge x_2 \ge\cdots\ge x_r \ge2$, for any $0<\epsilon<x_r$, we have
\begin{eqnarray}
&&\limsup_{m \to\infty} \frac{1}{m} \log Q_m \bigl(
\xi_{1}^{m,0} \ge x_1,\ldots,
\xi_{r}^{m,0} \ge x_r \bigr)
\nonumber
\\
&&\quad  \le\limsup_{m \to\infty} \frac{1}{m} \log
\bigl(Q_m \bigl(\xi_{1}^{m} \ge
x_1-\epsilon,\ldots,\xi_{r}^{m} \ge
x_r-\epsilon \bigr) +\mathbb{P} (g_m\ge\epsilon ) \bigr)
\nonumber
.
\end{eqnarray}
Moreover,\vspace*{1pt}
\[
Q_m \bigl(\xi_{1}^{m} \ge x_1-
\epsilon,\ldots,\xi_{r}^{m} \ge x_r-\epsilon
\bigr) = \exp\bigl\{ -m\bigl(I_r(x_1-
\epsilon,x_2-\epsilon,\ldots ,x_r-\epsilon)+\mathrm{o}(1)
\bigr)\bigr\},
\]
where $\mathrm{o}(1)$ goes to 0 as $m$ goes to infinity.
So for fixed $0<\epsilon<x_r$,\vspace*{1pt}
\[
\frac{\mathbb{P}(g_m\ge\epsilon)}{Q_m (\xi_{1}^{m} \ge x_1-\epsilon
,\ldots,\xi_{r}^{m}
\ge x_r-\epsilon )}\to0,\qquad  m\to\infty,
\]
hence,
%
%e5.11 #&#
\[
\limsup_{m \to\infty} \frac{1}{m} \log Q_m \bigl(
\xi_{1}^{m,0} \ge x_1,\ldots,
\xi_{r}^{m,0} \ge x_r \bigr)\le-I_r(x_1-
\epsilon ,x_2-\epsilon,\ldots,x_r-\epsilon).
\nonumber
\]
Likewise,
%
%e5.12 #&#
\[
\liminf_{m \to\infty} \frac{1}{m} \log Q_m \bigl(
\xi_{1}^{m,0} \ge x_1,\ldots,
\xi_{r}^{m,0} \ge x_r \bigr)
\ge-I_r(x_1+\epsilon,x_2+\epsilon,
\ldots,x_r+\epsilon),
\]
and letting $\epsilon\to0$, the continuity of the rate
function leads to \eqref{multi2:finite}.
\end{pf}

For any $\mu\in\mathcal{P}(\mathbb{R})$, construct a discrete
approximation via
%
%e5.13 #&#
\begin{equation}\label{traceless:discrete}
x_{i}^m=\inf \biggl\{x \in\mathbb{R}\dvtx \mu\bigl((-\infty,x]\bigr)
\ge\frac{i}{m+1} \biggr\}, \qquad 1\le i \le m,
\end{equation}
and $\mu^m=\sum_{i=1}^m \delta_{x_i^m}/m$ (note that the choice of the
length $1/(m+1)$ of the intervals rather that $1/m$ is only made in
order to
ensure that $x_m^m$ is finite).
Using these discrete constructions, set:
%
%e5.14 #&#
\begin{equation}
\label{A} \mathcal{X}= \Biggl\{\mu\in\mathcal{P}(\mathbb{R})\dvtx
\frac{1}{\sqrt{m}}\sum_{i=1}^{m}x_i^m
\to0, \mbox{ as } m \to\infty \Biggr\}
\end{equation}
and
%
%e5.15 #&#
\begin{equation}
\mathcal{P}_0(\mathbb{R})= \biggl\{\mu\in\mathcal{P}(\mathbb{R})\dvtx
\int x \mu(\mathrm{d}x)=0 \biggr\}.
\end{equation}

Since the condition in \eqref{A} ensures that
$\mu$ has mean zero, it is clear that $\mathcal{X}$ is a proper subset of
$\mathcal{P}_0(\mathbb{R})$. With the above, and the arguments and
results in \cite{BG},
the large deviation principle for the
spectral measure of the traceless GUE follow:
%
%th5.2 #&#
\begin{theorem}\label{T:traceless}
The spectral measure $\hat{\mu}_0^m=\sum_{i=1}^m\delta_{\xi_i^0}/m$
satisfies a large deviation principle on $\mathcal{X}$ in the scale
$m^2$ and
with the good rate function $I$.
\end{theorem}
\begin{pf} Since this proof closely follows \cite{BG}, it is
just sketched here. Write the density of the eigenvalues as:
\begin{eqnarray*}
&&Q_m\bigl(\mathrm{d}\xi_{1}^0,\mathrm{d}
\xi_{2}^0,\ldots, \mathrm{d}\xi_{m}^0
\bigr)
\nonumber
\\
&&\quad =\frac{\sqrt{2\uppi}}{Z_m}\exp \biggl(-m^2\iint_{x\ne y}f(x,y)
\hat{\mu }_0^m(\mathrm{d}x)\hat{\mu}_0^m(
\mathrm{d}y) \biggr)\prod_{i=1}^m
\mathrm{e}^{-\sfrac{{\xi
_i^0}^2}{2}}\,\mathrm{d}\xi_1^0\cdots\,\mathrm{d}
\xi_{m-1}^0,
\nonumber
\end{eqnarray*}
where $\xi_{m}^0=-\sum_{i=1}^{m-1}\xi_i^0$ and
\[
f(x,y)=\tfrac{1}{4}\bigl(x^2+y^2\bigr)-\log|x-y|.
\]

Let $\bar{Q}_m$ be the non-normalized positive measure $\bar
{Q}_m=Z_mQ_m/\sqrt{2\uppi}$.
Via Stirling's formula,
%
%e5.16 #&#
\begin{equation}
\label{traceless:factor} \lim_{m \to\infty}\frac{1}{m^2}\log
\frac{\sqrt{2\uppi}}{Z_m} = \frac{1}{2}-\int_0^1
x\log x \,\mathrm{d}x =\frac{3}{4},
\end{equation}
so if under $\bar{Q}_m$, $\hat{\mu}_0^m$ satisfies a large deviation
with rate function
%
%e5.17 #&#
\begin{equation}\label{traceless:newrate}
J(\mu)=\iint f(x,y)\mu(\mathrm{d}x)\mu(\mathrm{d}y),
\end{equation}
then combined with \eqref{traceless:factor}, this will lead to
the statement of the theorem.

First, observe that for any Borel subset $A\subset\mathcal{X}$, any
$N\in\mathbb{R}^+$,
%
%e5.18 #&#
\begin{equation}\label{traceless:upperbound}
\limsup_{m\to\infty}\frac{1}{m^2}\log \bigl(
\bar{Q}_m\bigl(\hat{\mu}_0^m \in A\bigr)
\bigr) \le-\inf_{\mu\in A} \biggl(\iint\bigl(f(x,y)\wedge N\bigr)
\mu(\mathrm{d}x)\mu(\mathrm{d}y) \biggr).
\end{equation}

Moreover, from arguments as in \cite{BG}, the sequence
$(\hat{\mu}_0^m)_{m\in\mathbb{N}}$ is exponentially tight
under $\bar{Q}_m$ on $\mathcal{X}$. So we just need to prove that
$(\hat{\mu}_0^m)_{m\in\mathbb{N}}$ satisfies a weak large deviation
principle with rate function $J(\mu)$ under the measure $\bar{Q}_m$.
The upper bound is clear. Indeed, $\mu\to\iint(f(x,y)\wedge N) \mu
(\mathrm{d}x)\mu(\mathrm{d}y)$ is
continuous on $\mathcal{X}$, therefore \eqref{traceless:upperbound}
implies that for any
$\mu\in\mathcal{X}$,
\[
\limsup_{\delta\to0}\limsup_{m \to\infty}
\frac{1}{m^2}\log \bigl(\bar {Q}_m\bigl(\hat{
\mu}_0^m \in B(\mu,\delta) \bigr) \bigr)\le-\iint
\bigl(f(x,y)\wedge N\bigr) \mu(\mathrm{d}x)\mu(\mathrm{d}y),
\]
where $B(\mu,\delta)$ is the open ball of center $\mu$ and radius
$\delta$,
with respect to the distance given by
\[
d(\mu_1,\mu_2)=\sup_{g\in \Lip_{\mathrm{b}}(1)}\biggl
\llvert \int g \,\mathrm{d}\mu_1-\int g \,\mathrm{d}\mu
_2\biggr\rrvert ,\qquad  \mu_1, \mu_2 \in
\mathcal{X},
\]
where for some fixed $b > 0$,
\[
\Lip_{\mathrm{b}}(1)=\bigl\{g\dvtx \mathbb{R}\to\mathbb{R}\dvtx \|g\|_{\Lip}
\le1, \|g\|_{\infty
}\le b\bigr\},
\]
are bounded Lipschitz functions. Then, by monotone convergence,\vspace*{-1pt}
%
%e5.19 #&#
\begin{equation}\label{traceless:upperbound1}
\limsup_{\delta\to0}\limsup_{m \to\infty}
\frac{1}{m^2}\log \bigl(\bar {Q}_m\bigl(\hat{
\mu}_0^m \in B(\mu,\delta) \bigr) \bigr)\le-\iint f(x,y)
\mu (\mathrm{d}x)\mu(\mathrm{d}y),
\end{equation}
finishing the proof of the upper bound.

To prove the lower bound, let $\nu\in\mathcal{X}$.
Since $I(\nu)=+\infty$ if $\nu$ has an atom, assume without loss of generality
that $\nu$ is atomless. As in \eqref{traceless:discrete}, let
$\nu^m=\sum_{i=1}^m \delta_{x_i^m}/m$. Now, as $m\to\infty$, $\nu^m$
converges weakly, with probability one,
towards $\nu$. Hence, for any $\delta>0$ and $m$ large enough, setting
$\Delta_m:=\{\xi_{1}^0 \le\xi_{2}^0 \le\cdots\le\xi_{m}^0\}$,
%
%e5.20 #&#
%e5.21 #&#
\begin{eqnarray}\label{traceless:probability}\label{E:low}
&&\hspace*{-35pt}\bar{Q}_m\bigl(\hat{\mu}_0^m \in B(\nu,
\delta)\bigr)\nonumber\\
&&\hspace*{-35pt}\quad \ge\bar{Q}_m \biggl( \biggl\{ \max_{1\le i \le m-1}\bigl|
\xi_i^0-x_i^m\bigr|<
\frac{\delta}{2\sqrt{m}} \biggr\} \cap\Delta_m \biggr)\nonumber
\\
&&\hspace*{-35pt}\quad \ge\int_{\mathcal{T}(\xi_1,\ldots,\xi_m)} \exp \Biggl(-\frac{m}{2}\sum
_{i=1}^m \bigl(\xi_i+x_i^m
\bigr)^2 \Biggr)\prod_{1\le i<j\le m}\bigl|
\xi_i-\xi_j+x_i^m-x_j^m\bigr|^2
\prod_{i=1}^{m-1}\,\mathrm{d}
\xi_i\nonumber
\\
&&\hspace*{-35pt}\quad \ge\prod_{i+1<j}\bigl|x_i^m-x_j^m\bigr|^2
\times\prod_{i=1}^{m-1}\bigl|x_{i+1}^m-x_i^m\bigr|
\exp \Biggl\{-\frac{m}{2}\sum_{i=1}^{m-1}
\biggl(\bigl|x_i^m\bigr|+\frac{\delta}{\sqrt{m}} \biggr)^2
\Biggr\}
\nonumber
\\
&&\hspace*{-35pt}\qquad {}\times\bigl|x_m^m-x_{m-1}^m\bigr|\exp
\Biggl( -m \Biggl(\sum_{i=1}^{m-1}x_i^m
\Biggr)^2 -m^2 \delta^2 \Biggr) \int
_{\mathcal{T}(\xi_1,\ldots,\xi_m)} \prod_{i=1}^{m-2}|
\xi_{i+1}-\xi_{i}| \prod_{i=1}^{m-1}
\,\mathrm{d}\xi_i,
\end{eqnarray}
where
\begin{eqnarray*}
&&\mathcal{T}(\xi_1,\ldots,\xi_m)
\nonumber
\\
&&\quad := \Biggl\{\max_{1\le i \le m-1}|\xi_i|<
\frac{\delta}{2\sqrt{m}}, \xi _{1} \le\xi_{2} \le\cdots\le
\xi_{m}, \sum_{i=1}^m
\xi_i+\sum_{i=1}^m
x_i^m=0 \Biggr\}.
\nonumber
\end{eqnarray*}
The last term on the right-hand side of \eqref{E:low} can be lower-bounded
by changing variables: $\xi_1=x_1$ and $\xi_{i}-\xi_{i-1}=x_{i}, 2\le i
\le m-1$.
Next, let
\begin{eqnarray*}
&&\mathcal{R}(x_1,\ldots,x_{m-1})
\nonumber
\\
&&\quad := \biggl\{-\frac{\delta}{2\sqrt{m}}\le x_1 \le-\frac{\delta}{4\sqrt{m}}, 0
\le x_i \le-\frac{\delta}{4m^2}, 2 \le i \le m-1 \biggr\}.
\nonumber
\end{eqnarray*}
Recalling that, $\sum_{i=1}^{m}x_i^m/m \to0$,
%
%e5.22 #&#
\begin{eqnarray}
\int_{\mathcal{T}(\xi_1,\ldots,\xi_m)}\prod_{i=1}^{m-2}|
\xi_{i+1}-\xi _{i}|\prod_{i=1}^{m-1}
\,\mathrm{d}\xi_i&\ge&\int_{\mathcal{R}(x_1,\ldots
,x_{m-1})}\prod
_{i=2}^{m-1}|x_i|\prod
_{i=1}^{m-1}\,\mathrm{d}x_i
\nonumber
\\[-8pt]\\[-8pt]
&\ge&\frac{\delta}{4\sqrt{m}} \biggl(\frac{1}{2} \biggl(\frac{\delta
}{4m^2}
\biggr)^2 \biggr)^{m-2}.\nonumber
\end{eqnarray}
Hence,
%
%e5.23 #&#
%e5.24 #&#
\begin{eqnarray}
\label{traceless:lowerbound1}
&&\bar{Q}_m\bigl(\hat{\mu}_0^m
\in B(\nu,\delta)\bigr) \nonumber\\
&&\quad \ge\prod_{i+1<j}\bigl|x_i^m-x_j^m\bigr|^2
\prod_{i=1}^{m-1}\bigl|x_{i+1}^m-x_i^m\bigr|
\exp \Biggl(-\frac{m}{2}\sum_{i=1}^{m}
\bigl(x_i^m\bigr)^2 \Biggr)
\\
&&\qquad {}\times\bigl|x_m^m-x_{m-1}^m\bigr|
\frac{\delta}{4\sqrt{m}} \biggl(\frac{1}{2} \biggl(\frac{\delta}{4m^2}
\biggr)^2 \biggr)^{m-2} \exp \Biggl(-\sqrt{m}\delta \sum
_{i=1}^m\bigl|x_i^m\bigr|-
\delta^2 \Biggr).\nonumber
\end{eqnarray}
Now by arguments as in \cite{BG},
%
%e5.25 #&#
\begin{equation}\label{traceless:lowerbound5}
\liminf_{\delta\to0}\liminf_{m \to\infty}
\frac{1}{m^2}\log \bigl(\bar {Q}_m\bigl(\hat{
\mu}_0^m \in B(\nu,\delta) \bigr) \bigr)\ge-\iint f(x,y)
\nu (\mathrm{d}x)\nu(\mathrm{d}y).
\end{equation}
Combining \eqref{traceless:upperbound1} and \eqref
{traceless:lowerbound5}, establishes
the weak large deviation principle, finishing the proof of the theorem.
\end{pf}

We are now ready to give the large deviations for $\xi_{1}^{m,0}$ when
on the left
of its mean. To do so, let us introduce some notations:
Let $\mathcal{M}((-\infty,x])$ be the set of all probability measures
on $(-\infty,x],x\in\mathbb{R}$,
let $\mathcal{M}_{\mathcal{X}}((-\infty,x])=\mathcal{M}((-\infty,x])\cap
\mathcal{X}$, and let\vspace*{2pt}
$\mathcal{M}_0((-\infty,x])=\mathcal{M}((-\infty,x])\cap\mathcal
{P}_0(\mathbb{R})$.
Since $\{\xi_{1}^{m,0}\le x\}=\{\hat{\mu}_0^m \in\mathcal{M}_{\mathcal
{X}}((-\infty,x]) \}$,
then for any $x\le2$,
%
%e5.26 #&#
\begin{equation}
\label{E:coro1} \lim_{m \to\infty} \frac{1}{m^2} \log\mathbb{P}
\bigl(\xi_{1}^{m,0} \leq x \bigr)=-\inf_{\mu\in\mathcal{M}_{\mathcal{X}}((-\infty,x])}I(
\mu).
\end{equation}
For each $x\in\mathbb{R}$, let
%
%e5.27 #&#
\begin{equation}
\label{E:appendix} K(x)=\inf_{\mu\in\mathcal{M}_{0}((-\infty,x])}I(\mu).
\end{equation}
When $x\ge2$, the semicircular law $\sigma$ is both in $\mathcal
{M}_{\mathcal{X}}((-\infty,x])$ and $\mathcal{M}_{0}((-\infty,x])$, and
so $\inf_{\mu\in\mathcal{M}_{\mathcal{X}}((-\infty,x])}I(\mu
)=K(x)=I(\sigma)=0$.
Moreover, when $x \le0$, and since both $\mathcal{M}_{\mathcal
{X}}((-\infty,\allowbreak x])$
and $\mathcal{M}_{0}((-\infty,x])$ are empty,
$\inf_{\mu\in\mathcal{M}_{\mathcal{X}}((-\infty,x])}I(\mu)=K(x)=I(\sigma
)=+\infty$.

When $0 < x \le2$, and from arguments as in \cite{Ib}, it is next
shown that $K$ is continuous. Indeed, for any $y<0$ and $0 < x \le2$, let
%
%e5.28 #&#
\begin{equation}
J_{\mu}(y, x)=\frac{1}{2}\int_{y}^{x}
u^2 \mu(\mathrm{d}u) -\int_{y}^{x}
\int_{y}^{x}\log|u-t|\mu(\mathrm{d}u)\mu(
\mathrm{d}t)-\frac{3}{4} ,
\end{equation}
and let $\nu_x$ be the minimizer of $I(\mu)$ on $\mathcal
{M}_{0}((-\infty,x])$.
Then, for any $0 < \epsilon< x$,
%
%e5.29 #&#
\begin{equation}
\label{A:leftconti} K(x)\le K(x-\epsilon)\le\frac{J_{\nu_x}(y_\epsilon,x-\epsilon)}{\nu
_x^2([y_\epsilon,x-\epsilon])},
\end{equation}
where $y_\epsilon$ is the value for which
\[
\int_{y_\epsilon}^{x-\epsilon} t \,\mathrm{d}\nu_x
(t)=0.
\]
Since, as $\epsilon\to0$, the right-hand side of \eqref{A:leftconti}
converges to
$K(x)$, $K$ is left continuous.

To show the right continuity, note that by a simple change of variables,
\[
K(x)=\inf_{\mu\in\mathcal{M}_{0}((-\infty,x+\epsilon])}J_{\mu}^{\epsilon}(x),
\]
where
\[
J_{\mu}^{\epsilon}(x)=\frac{1}{2}\int_{-\infty}^{x+\epsilon}
(u-\epsilon )^2 \mu(\mathrm{d}u)-\int_{-\infty}^{x+\epsilon}
\int_{-\infty}^{x+\epsilon}\log |u-t|\mu(\mathrm{d}u)\mu(
\mathrm{d}t)-\frac{3}{4}.
\]
Therefore,
\[
0\le K(x)-K(x+\epsilon)\le J_{\nu_{x+\epsilon}}^{\epsilon
}(x)-K(x+\epsilon)=
\frac{\epsilon^2}{2},
\]
and the right continuity of $K$ follows.
Likewise, $\inf_{\mu\in\mathcal{M}_{\mathcal{X}}((-\infty,x])}I(\mu)$ is
right-continuous with respect to $x$.

Next we need a result which, when combined with \eqref{E:coro1}, gives
%
%e5.30 #&#
\begin{equation}
\label{E:coro22} \lim_{m \to\infty} \frac{1}{m^2} \log\mathbb{P}
\bigl(\xi_{1}^{m,0} \leq x \bigr)=-K(x),
\end{equation}
for any $x\le2$. This is the purpose of our next lemma whose statement
as well as
proof benefited from Ionel Popescu help.
%
%le5.1 #&#
\begin{lemma}\label{lemma_appendix}
For any $x \in\mathbb{R}$,
%
%e5.31 #&#
\begin{equation}
\label{lemmaeq} \inf_{\mu\in\mathcal{M}_{\mathcal{X}}((-\infty,x])}I(\mu)=K(x).
\end{equation}
\end{lemma}
\begin{pf} For $x \ge2$, both sides of \eqref{lemmaeq} are equal
to zero and so
we just need to consider the case $x <2$.
First, since $\mathcal{X}$ is a proper subset of
$\mathcal{P}_0(\mathbb{R})$,
%
%e5.32 #&#
\begin{equation}
\label{E:lemma_1} K(x)\le\inf_{\mu\in\mathcal{M}_{\mathcal{X}}((-\infty,x])}I(\mu).
\end{equation}
Next, let us show that
%
%e5.33 #&#
\begin{equation}
\label{E:lemma_2} K(x)\ge\inf_{\mu\in\mathcal{M}_{\mathcal{X}}((-\infty,x])}I(\mu).
\end{equation}
By Theorem~1.10 and Theorem~1.11 of Chapter IV of \cite{SF},
there exists a unique probability measure, $\mu_0$,
minimizing $I(\mu)$, for all $\mu\in\mathcal{M}_{0}((-\infty,x])$, and its
support is an interval, $[a,b]$, with $b \le x $.
Since $\mu_0$ is atomless, its distribution function $F$ is continuous,
increasing
with $F(a)=0$ and $F(b)=1$. Moreover, since $\mu_0$ has mean zero,
$\int_0^1 F^{-1}(x) \,\mathrm{d}x=0$, where $F^{-1}$, the inverse of $F$, is continuous
and increasing on $[0,1]$, with $F^{-1}(0)=a$ and $F^{-1}(1)=b$.

Now for any integer $n \ge2$, construct an approximation to $F^{-1}$
as follows:
For $i/n \le x \le(i+1)/n$, let
\[
G_{n}^{+}(x)= %
\cases{ n \biggl(F^{-1}
\biggl(\displaystyle \frac{i+2}{n}\biggr) \biggl(x-\displaystyle \frac{i}{n}
\biggr)+F^{-1}\biggl(\displaystyle \frac{i+1}{n}\biggr) \biggl(
\displaystyle \frac{i+1}{n}-x \biggr) \biggr), &\quad \mbox{if $0\le i \le n-2$,}
\cr
b+x-
\displaystyle \frac{i}{n}, &\quad \mbox{if $i=n-1$,} } %
\]
and let,
\[
G_{n}^{-}(x)= %
\cases{ n \biggl(F^{-1}
\biggl(\displaystyle \frac{i}{n}\biggr) \biggl(x-\displaystyle \frac{i}{n}
\biggr)+F^{-1}\biggl(\displaystyle \frac{i-1}{n}\biggr) \biggl(
\displaystyle \frac{i+1}{n}-x \biggr) \biggr), &\quad \mbox{if $1\le i \le n-1$,}
\cr
a+x-
\displaystyle \frac{i+1}{n}, &\quad \mbox{if $i=0$.} } %
\]
From this construction, $\int_0^1 G_{n}^{+}(x) \,\mathrm{d}x >0$ and $\int_0^1
G_{n}^{-}(x) \,\mathrm{d}x <0$.
Next, let
\[
\gamma_n^{+}=\frac{-\int_0^1 G_{n}^{-}(x) \,\mathrm{d}x}{\int_0^1 G_{n}^{+}(x) \,\mathrm{d}x
-\int_0^1 G_{n}^{-}(x) \,\mathrm{d}x}, \qquad \gamma_n^{-}
=\frac{\int_0^1 G_{n}^{+}(x) \,\mathrm{d}x}{\int_0^1 G_{n}^{+}(x) \,\mathrm{d}x-\int_0^1
G_{n}^{-}(x) \,\mathrm{d}x},
\]
and let
\[
G_n(x)=\gamma_n^{+}G_{n}^{+}(x)+
\gamma_n^{-}G_{n}^{-}(x).
\]
Then,
\[
\int_0^1 G_n(x) \,\mathrm{d}x=0,
\]
and since $G_n$ is piece-wise linear, it is Lipschitz.
Let $\mu_n$ be the probability measure whose distribution function is
$G_n^{-1}$.
The Lipschitz continuity of $G_n$ yields that $\mu_n \in\mathcal{X}$,
for any $n \ge2$.
From its very construction, $\mu_n$ is supported on $[a-1/n,b+1/n]$,
and $\mu_n$ converges to
$\mu_0$ weakly, as $n\to\infty$, and thus
%
%e5.34 #&#
\begin{equation}
\label{E:appen-rate1} \lim_{n \to\infty}\int x^2
\mu_n(\mathrm{d}x) = \int x^2 \mu_0 (
\mathrm{d}x).
\end{equation}
For the second term on the right-hand side of \eqref{traceless:ratefunction},
%
%e5.35 #&#
\begin{equation}
\iint\log|x-y|\mu(\mathrm{d}x)\mu(\mathrm{d}y)=2 \iint_{x < y} \log(y-x)
\mu(\mathrm{d}x)\mu(\mathrm{d}y),
\end{equation}
let
%
%e5.36 #&#
\begin{equation}
\frac{1}{n^2} \sum_{i<j} \log \biggl(
F^{-1} \biggl(\frac{j+1}{n} \biggr)-F^{-1} \biggl(
\frac{i}{n} \biggr) \biggr)+\frac{1}{2n^2} \sum
_{i=0}^{n-1} \log \biggl( F^{-1} \biggl(
\frac{i+1}{n} \biggr)-F^{-1} \biggl(\frac{i}{n} \biggr)
\biggr)
\end{equation}
and
%
%e5.37 #&#
\begin{equation}
\frac{1}{n^2} \sum_{i<j}\log \biggl(
G_n \biggl(\frac{j+1}{n} \biggr)-G_n \biggl(
\frac{i}{n} \biggr) \biggr)+\frac{1}{2n^2}\sum
_{i=0}^{n-1}\log \biggl( G_n \biggl(
\frac{i+1}{n} \biggr)-G_n \biggl(\frac{i}{n} \biggr)
\biggr),
\end{equation}
be respectively Riemann sums approximations of $\iint_{x < y} \log(y-x)
\mu_0(\mathrm{d}x)\mu_0(\mathrm{d}y)$
and\linebreak[4]  $\iint_{x < y} \log(y-x) \mu_n(\mathrm{d}x)\mu_n(\mathrm{d}y)$. For any $0 \le i\le j
\le n-1$,
%
%e5.38 #&#
\begin{eqnarray}
&&\log \biggl(G_n \biggl(\frac{j+1}{n} \biggr)-G_n
\biggl(\frac{i}{n} \biggr) \biggr)
\nonumber
\\[-8pt]\\[-8pt]
&&\quad \ge\gamma_n^{+}\log \biggl(G_{n}^{+}
\biggl(\frac{j+1}{n} \biggr)-G_{n}^{+} \biggl(
\frac{i}{n} \biggr) \biggr) +\gamma_n^{-}\log
\biggl(G_{n}^{-} \biggl(\frac{j+1}{n} \biggr)
-G_{n}^{-} \biggl(\frac{i}{n} \biggr) \biggr),\nonumber
\end{eqnarray}
and moreover, for any $1\le i \le j \le n-2$,
%
%e5.39 #&#
\begin{eqnarray}
&&\log \biggl(G_n \biggl(\frac{j+1}{n} \biggr)-
G_n \biggl(\frac{i}{n} \biggr) \biggr)
\nonumber\\
&&\quad \ge\gamma_n^{+}\log \biggl( F^{-1} \biggl(
\frac{j+2}{n} \biggr)-F^{-1} \biggl(\frac{i+1}{n} \biggr)\biggr)
\\
&&\quad \quad + \,\gamma_n^{-}\log \biggl( F^{-1} \biggl(
\frac{j}{n} \biggr)-F^{-1} \biggl(\frac{i-1}{n} \biggr)
\biggr).\nonumber
\end{eqnarray}

If $\iint_{x < y} \log(y-x) \mu_0(\mathrm{d}x)\mu_0(\mathrm{d}y)=-\infty$, \eqref
{E:lemma_2} is trivially true, so
let us assume that this last integral is finite. Moreover, since $\gamma
_n^{+}+\gamma_n^{-}=1$,
%
%e5.40 #&#
\begin{equation}
\label{E:appen-rate2} \liminf_{n \to\infty} \biggl({-\!\!}\iint\log|x-y|
\mu_n(\mathrm{d}x)\mu_n(\mathrm{d}y) \biggr)\le {-\!\!}\iint
\log|x-y|\mu_0(\mathrm{d}x)\mu_0(\mathrm{d}y),
\end{equation}
and combining \eqref{E:appen-rate1} and \eqref{E:appen-rate2},
\[
\liminf_{n \to\infty}I(\mu_n)\le I(\mu_0).
\]
Since $\mu_n$ is supported on $[a-1/n,b+1/n]$ and from the right
continuity (in $x$)
of $\inf_{\mu\in\mathcal{M}_{\mathcal{X}}((-\infty,x])}I(\mu)$,
\[
K(x)\ge\inf_{\mu\in\mathcal{M}_{\mathcal{X}}((-\infty,x])}I(\mu),
\]
which finishes the proof.
\end{pf}

To finish this Appendix, the large deviations for the first $r$
eigenvalues of the traceless GUE, when at least one of them is on the
left of the asymptotic mean,
is established:
%
%co5.1 #&#
\begin{corollary}\label{coro2}
For $x_r\le x_{r-1} \le\cdots\le x_1 $, and $x_r\le2$,
\[
\lim_{m \to\infty} \frac{1}{m^2} \log\mathbb{P} \bigl(
\xi_{1}^{m,0} \leq x_1,\ldots,
\xi_{r}^{m,0} \leq x_r \bigr)=-K(x_r).
\]
\end{corollary}
\begin{pf} Next, let $(\xi_{1}^{m,0},\xi_{2}^{m,0},\ldots, \xi
_{m}^{m,0})$ be the non-increasing rearrangement of $(\xi_{1}^{0},\xi
_{2}^{0},\ldots, \xi_{m}^{0})$, and set
\begin{eqnarray*}
L&:=&\mathbb{P} \bigl(\xi_{1}^{m,0} \leq x_1,
\ldots,\xi_{r}^{m,0} \leq x_r \bigr),
\nonumber
\\
M&:=&\mathbb{P} \bigl(\xi_{1}^{0} \leq x_1,
\ldots,\xi_{r}^{0} \leq x_r,\xi
_{r+1}^{0} \leq x_r,\ldots,\xi_{m}^{0}
\leq x_r \bigr).
\nonumber
\end{eqnarray*}
Then,
%
%e5.41 #&#
\begin{equation}\label{traceless:corollary}
M \le L\le\frac{m!}{(m-r+1)!(r-1)!}B \le m^r M,
\end{equation}
and therefore,
%
%e5.42 #&#
\begin{equation}\label{traceless:corollary1}
\lim_{m \to\infty} \frac{1}{m^2} \log L=\lim
_{m \to\infty} \frac
{1}{m^2} \log M.
\end{equation}
Changing variables:
\begin{eqnarray*}
\xi_i^0-(x_i-x_r)&=&
\eta_i, \qquad \mbox{for } 1\le i \le r-1,
\\
\xi_i^0&=&\eta_i, \qquad \mbox{for } r \le i \le
m,
\end{eqnarray*}
and so,
\[
M=\mathbb{P} (\eta_i \leq x_r, 1\le i \le m ).
\]
Considering the two measures $\sum_{i=1}^m\delta_{\xi_i^0}/m$ and
$\sum_{i=1}^m \delta_{\eta_i}/m$, for any bounded Lipschitz function
$g$ (with $\|g\|_{\Lip}\le1$), then as $m\to\infty$,
\[
\frac{1}{m}\Biggl\llvert \sum_{i=1}^m
g\bigl(\xi_i^0\bigr)-\sum_{i=1}^m
g(\eta_i)\Biggr\rrvert \le\frac{1}{m}\sum
_{i=1}^m\bigl|\xi_i^0-
\eta_i\bigr|\longrightarrow0.
\]
Therefore, $\sum_{i=1}^m\delta_{\xi_i^0}/m$ and
$\sum_{i=1}^m \delta_{\eta_i}/m$ are exponentially equivalent, and
Theorem~\ref{T:traceless} also applies to the latter (see Theorem~4.2.13 in \cite{DemboZ}). So from \eqref{E:coro22}, it follows that
\[
\lim_{m \to\infty} \frac{1}{m^2} \log M=-K(x_r),
\]
and \eqref{traceless:corollary1} finishes the proof.
\end{pf}

\end{appendix}

\section*{Acknowledgements}
Many thanks to Satya Majumdar
for suggesting that the methods of \cite{NMV} should
provide a closed form expression for the rate function obtained in
Theorem~\ref{P:tableautraceless}, to Ionel Popescu for his help with
Lemma~\ref{lemma_appendix} and its proof, and to a referee for detailed
comments.
Christian Houdr\'e's research was supported in part by the Simons
Foundation Grant, \#246283.

%\begin{supplement}%[id=suppA]
%\slink[doi]{10.3150/00-BEJXXXXSUPP} %[doi,text={...}] - jei reikia
%suskaldyti doi
%\sfilename{BEJ000\_supp.pdf}
%\end{supplement}

% imsref loaded by jurgita.kaciuliene, 2014-04-18 11:37:34

\printhistory

\end{document}